\documentclass[10pt,a4paper]{article}
\usepackage{t1enc}
\usepackage[utf8]{inputenc}
\usepackage{graphicx}

\usepackage{amsmath}
\usepackage{amssymb}
\usepackage{amsthm}
\usepackage{mathrsfs}
\usepackage[all]{xy}
\usepackage{enumerate}

\theoremstyle{plain}
\newtheorem{thm}{Theorem}[section]
\newtheorem{prop}[thm]{Proposition}

\newtheorem{cor}[thm]{Corollary}
\newtheorem{lem}[thm]{Lemma}

\theoremstyle{definition}
\newtheorem{dfn}{Definition}[section]

\theoremstyle{remark}
\newtheorem*{rem}{Remark}

\DeclareMathOperator{\Aut}{Aut}

\DeclareMathOperator{\Gal}{Gal}

\DeclareMathOperator{\Spec}{Spec}
\DeclareMathOperator{\Spf}{Spf}

\DeclareMathOperator{\ga}{\pi_1^{alg}}
\DeclareMathOperator{\gtop}{\pi_1^{top}}

\DeclareMathOperator{\gtemp}{\pi_1^{temp}}
\DeclareMathOperator{\glog}{\pi_1^{log}}
\DeclareMathOperator{\ggeom}{\pi_1^{log-geom}}
\DeclareMathOperator{\gtempgeom}{\pi_1^{temp-geom}}

\DeclareMathOperator{\Norm}{Norm}

\DeclareMathOperator{\KCov}{KCov}
\DeclareMathOperator{\KCovgeom}{KCov_{geom}}

\DeclareMathOperator{\GalKCov}{GalKCov}
\newcommand{\LGalKCov}{\mbb L\text{-}\GalKCov}
\DeclareMathOperator{\Covalg}{Cov^{alg}}
\DeclareMathOperator{\Covtop}{Cov^{top}}
\DeclareMathOperator{\Covtemp}{Cov^{temp}}
\DeclareMathOperator{\ket}{k\acute{e}t}

\DeclareMathOperator{\Btemp}{\mathcal B^{temp}}
\DeclareMathOperator{\BtempL}{\mathcal B^{temp,\mbb L}}

\DeclareMathOperator{\rk}{rk}
\DeclareMathOperator{\gp}{gp}
\DeclareMathOperator{\Ker}{Ker}
\DeclareMathOperator{\Coker}{Coker}
\DeclareMathOperator{\red}{red}
\DeclareMathOperator{\Str}{Str}

\DeclareMathOperator{\DDtemp}{DD_{temp}}
\DeclareMathOperator{\Dtop}{\mcal D_{top}}
\DeclareMathOperator{\Dtopan}{\mcal D_{top}^{an}}
\DeclareMathOperator{\Dtops}{\mcal D_{top}^{sp}}
\DeclareMathOperator{\Dtopgeom}{\mcal D_{top-geom}}

\DeclareMathOperator{\Ob}{Ob}
\DeclareMathOperator{\Hom}{Hom}

\DeclareMathOperator{\et}{\acute{e}t}

\DeclareMathOperator{\Zar}{Zar}

\DeclareMathOperator{\Ens}{Set}
\DeclareMathOperator{\Set}{Set}
\DeclareMathOperator{\tSet}{\text{-}Set}
\DeclareMathOperator{\tEns}{\text{-}Set}

\DeclareMathOperator{\fSet}{fSet}
\DeclareMathOperator{\Poset}{Poset }

\DeclareMathOperator{\Stab}{Stab}

\DeclareMathOperator{\an}{an}
\DeclareMathOperator{\Image}{Im}
\DeclareMathOperator{\codim}{codim}
\DeclareMathOperator{\Irr}{Irr}
\DeclareMathOperator{\sq}{\square}
\DeclareMathOperator{\sm}{sm}
\DeclareMathOperator{\tps}{tps}
\DeclareMathOperator{\op}{op}
\DeclareMathOperator{\Lim}{Lim}

\DeclareMathOperator{\injLim}{\underset{\longrightarrow}{\Lim}}

\DeclareMathOperator{\rklog}{rk^{log}}
\DeclareMathOperator{\Ke}{\mcal{K}e}
\newcommand{\mring}{\mathring}
\DeclareMathOperator{\triv}{tr}
\newcommand{\tr}{\triv}

\DeclareMathOperator{\Gm}{\mbf G_m}
\DeclareMathOperator{\fGm}{\fk G_m}
\DeclareMathOperator{\OutGptop}{OutGp_{top}}
\DeclareMathOperator{\Pt}{Pt}
\DeclareMathOperator{\nd}{nd}

\newcommand{\ie}{\emph{i.e.} }
\newcommand{\findem}{\end{proof}}
\newcommand{\dem}{\begin{proof}}

\newcommand{\da}{\begin{displaystyle}}
\newcommand{\db}{\end{displaystyle}}
\newcommand{\dar}{\downarrow}
\newcommand{\uar}{\uparrow}

\newcommand{\mcal}{\mathcal}
\newcommand{\mbf}{\mathbf}

\newcommand{\mbb}{\mathbb}
\newcommand{\fk}{\mathfrak}

\DeclareMathOperator{\C}{C}
\DeclareMathOperator{\Cgeom}{C_{geom}}
\DeclareMathOperator{\Ctop}{C_{top}}

\begin{document}

\title{Coverings in $p$-adic analytic geometry and log coverings II:\\
Cospecialization of the $(p')$-tempered fundamental group in higher dimensions}
\author{Emmanuel Lepage}
\maketitle
\tableofcontents

\section*{Introduction}
\subsection*{}
For a proper morphism of schemes $f: X\to Y$ with geometrically connected
fibers and a specialization $ 
\bar y_1 \to \bar y_2$ of geometric points of $Y$, A. Grothendieck  
has defined algebraic fundamental groups $\pi_1^{alg}(X_{\bar y_i})$  
and a specialization homomorphism  $\pi_1^{alg}(X_{\bar y_1})\to  
\pi_1^{alg}(X_{\bar y_2})$. Grothendieck's specialization theorem 
tells that this homomorphism is surjective if $f$ is separable and induces an  
isomorphism between the prime-to-$p$ quotients if $f$ is smooth  
(here, $p$ denotes the characteristic of $\bar y_2$, eventually $0$), cf. \cite[cor. X.2.4,
cor. X.3.9]{sga}.\\
In complex analytic geometry, a smooth and proper morphism is locally a trivial
fibration of real differential manifolds, so that, in particular, all the
fibers are homeomorphic, and thus have isomorphic (topological) fundamental
groups.\\

The aim of this paper is to find some analogous result to the specialization
theorem of Grothendieck in the frame of non archimedean analytic
geometry. We already studied the one dimensional case in dimension 1
in~\cite{cosp1}. The results obtained here require some assumptions that did
not appear in~\cite{cosp1}, so that the results of~\cite{cosp1} will not be
covered by the results of this article.\\
What will play the role of the fundamental group in the non archimedean
setting will be the \emph{tempered fundamental group} of Y. André.\\
The tempered fundamental group was introduced in~\cite[part III]{andre1} as
a sort of non archimedean analog of the topological fundamental group of complex algebraic
varieties; its profinite completion coincides with Grothendieck's algebraic
fundamental group; it has itself many infinite discrete quotients in
general.\\
Since the analytification (in the sense of V. Berkovich or of rigid geometry)
of a finite étale covering of a $p$-adic algebraic variety is not necessarily a topological covering, André had to use a slightly larger notion of covering.
He defined tempered coverings, which are \'etale coverings in
the sense of A.J. de Jong (that is to say that locally on the Berkovich
topology, it is a direct sum of finite coverings) such that, after pulling
back by some finite \'etale covering, they become topological coverings (for the
Berkovich topology). The \emph{tempered fundamental group} is the prodiscrete group that classifies those tempered coverings.  To give a more handful description, if
one has a sequence of pointed finite Galois connected coverings $((S_i,s_i))_{i\in \mbf N}$ such that the corresponding pointed pro-covering of $(X,x)$ is the
universal pro-covering of $(X,x)$, and if $(S^{\infty}_i,s^{\infty}_i)$ is a
universal topological covering of $S_i$, the tempered fundamental group of
$X$ can be seen as $\gtemp(X,x)=\varprojlim_i \Gal(S^{\infty}_i/X)$.
Therefore, to understand the tempered fundamental group of a variety, one
has to understand the topological behavior of its finite \'etale
coverings.\\

One problem which appears at once in looking for some non archimedean analog of
Grothendieck's specialization theorems is that there are in general no non trivial
specializations between distinct points of a non archimedean analytic
(Berkovich or rigid) space: for example a separated Berkovich space has a
Hausdorff underlying topological space, so that if there is a
cospecialization (for the Berkovich topology, the étale topology\dots)
between two geometric points of a Berkovich space, the two geometric points
must have the same underlying point. Thus we will study specialization
through a reduction (with good enough properties).\\
 
We want to understand how the tempered fundamental group of the geometric fibers of a
smooth family varies. Let us for instance consider a family of elliptic curves. The tempered
fundamental group of an elliptic curve over a complete algebraically closed
non archimedean closed field is $\widehat{\mbf Z}^2$ if it has
good reduction, and $\widehat{\mbf Z}\times \mbf Z$ if it is a Tate
curve. In particular, by looking at a moduli space of stable pointed elliptic curves
with level structure\footnote{to avoid stacks. However, the cospecialization
  homomorphisms we will construct will be local for the étale topology of
  the special fiber of the base. Thus, the fact that the base is a
  Deligne-Mumford stack is not really a problem.}, the tempered fundamental group (or any reasonable
$(p')$-version) cannot be constant.\\
Moreover, if one looks at the moduli space over $\mbf Z_p$, and one
considers a curve $E_1$ with bad reduction and a curve $E_0$ with generic
reduction (hence good reduction), there cannot be a morphism
$\gtemp(E_0)\to\gtemp(E_1)$ which induces Grothendieck's specialization on
the profinite completion, although the reduction point corresponding to
$E_1$ specializes to the reduction point corresponding to $E_0$. Therefore
there cannot be any reasonable specialization theory.\\
On the contrary, if
one has two geometric points $\eta_1$ and $\eta_2$ of the moduli space such
that the reduction of $\eta_1$ specializes to the reduction of $\eta_2$,
then $E_{\eta_1}$ has necessarily better reduction than $E_{\eta_2}$ and
there is some morphism
$\gtemp(E_{\eta_2})\to\gtemp(E_{\eta_1})$ that induces an isomorphism on
the profinite completions. Thus we want to look for a \emph{cospecialization} of the tempered fundamental group.\\
 
The topological behavior of a general finite étale covering is too
pathological to hope to have a simple cospecialization theory without
adding a $(p')$ condition on the coverings: for example
two Mumford curves over some finite extension of $\mbf Q_p$ with isomorphic
geometric tempered fundamental group have the same metrized graph of the stable reduction. Thus even if two Mumford curves
have isomorphic stable reduction (and thus the point corresponding to their stable reduction is the same), they do not have isomorphic
tempered fundamental group in general. Thus we will only
study here finite coverings that are dominated by a finite Galois covering
whose order is prime to $p$, where $p$ is the residual characteristic
(which can be $0$; such
a covering will be called a $(p')$-finite covering). Then, it becomes
natural to introduce a $(p')$-tempered fundamental group which classifies
tempered coverings that become topological coverings after pullback along
some $(p')$-finite covering. It should be remarked that this
$(p')$-tempered fundamental group cannot be in general recovered from the
tempered fundamental group.\\
The $(p')$-tempered fundamental group of a curve was already studied by
S. Mochizuki in~\cite{mochi}. It can be described in terms of a graph of profinite
groups. From this description, one easily sees that the isomorphism class
of the $(p')$-tempered
fundamental group of a $p$-adic curve depends only of the stratum of the
Knudsen stratification of the moduli space of stable curves in which the
stable reduction lies. Moreover if one has two strata $x_1$ and $x_2$ in the moduli space
of stable curves such that $x_1$ is in the closure of $x_2$, one can easily
construct  morphisms from the graph of groups
corresponding to $x_1$ to the graph of groups corresponding to $x_2$
(inducing morphisms of tempered fundamental groups which induce 
isomorphisms of the pro-$(p')$ completions).\\

\subsection*{}
We shall study the following general situation. Let $K$ be a discretely
valued field, $O_K$ be its valuation ring, $k$ be its residue field and $p$
its characteristics (which can be $0$).
Let $X\to Y$ be a pluristable (for example
semistable) morphism of schemes over $O_K$, which is smooth over $K$.\\
If $\eta_1$ is a (Berkovich) point of the
generic fiber of $Y$, we first want to describe the geometric
$(p')$-tempered fundamental group of $Y_{\eta_1}$ in terms of $Y_{s_1}$
where $s_1$ is the reduction of $\eta_1$ (to be sure that this reduction
exists we thus have to assume $Y$ to be proper over $O_K$; otherwise, we
have to consider only points $\eta_1$ in the tube of the
special fiber of $Y$). Let us make sure at first that we can get such a
description for the pro-$(p')$ completion, \emph{i.e.} the algebraic
fundamental group. One cannot apply directly Grothendieck's specialization
theorems since the special fiber is not smooth but only
pluristable. Indeed, a pro-$(p')$ geometric covering of the generic fiber will in
generally only induce a Kummer covering on the special fiber. These are
more naturally described in terms of \emph{log geometry}. To do so we will have to
assume the field $\mcal H(\eta_1)$ to be with discrete valuation in order
to get log schemes with good finiteness properties (more precisely to be
fs). Then, one can endow $X_{s_1}$
with a natural log structure such that the pro-$(p')$ fundamental group of
$X_{\eta_1}$ is isomorphic to a pro-$(p')$ log fundamental group (as defined
in~\cite{ill}). To try to describe the $(p')$-tempered fundamental group,
one has to describe the topological behavior of any $(p')$-algebraic
covering of $X_{\eta_1}$. Berkovich, in \cite{berk2}, constructed a combinatorial object
(more precisely a \emph{polysimplicial set}) depending only on $X_{s_1}$,
such that the Berkovich generic fiber $X_{\eta_1}$ is naturally
homotopically equivalent to the geometric realization of this combinatorial
object, thus generalizing the case of curves with semistable reduction,
where the homotopy type of the generic fiber can be naturally described in terms of
the graph of this semistable reduction. We will extend such a description to
our log coverings. This will enable us to define a
$(p')$-tempered fundamental group of our log reduction, which is isomorphic
to the tempered fundamental group of the generic fiber. In particular:
\begin{thm}The $(p')$-tempered fundamental group of $X_{\eta_1}$ only
  depends on the log reduction $X_{s_1}$.\end{thm}

Once we have a definition for the log geometric tempered fundamental group
$\gtempgeom(X_{s_1})$ of the log
fibers in the special locus of $Y$, one can reformulate our cospecialization
problem only in terms of this special locus.\\
We will now assume $Y$ to be log smooth over $O_K$ (this endows $Y$ with a
canonical stratification~; we did not used such an assumption for the case
of curves in~\cite{cosp1}).\\
We will also have to make an assumption on the combinatorial behavior of
the geometric fibers of $X\to Y$. More precisely, the polysimplicial set associated
with those geometric fibers will be assumed to be interiorly free. This is for example
the case if $X\to Y$ is strictly polystable or if $X\to Y$ is of relative
dimension 1 (which explains why such a condition did not appear in~\cite{cosp1}). We will prove the following:
\begin{thm}\label{sp} Let $\eta_1$ and $\eta_2$ be two Berkovich points
  with discrete valuation fields of
  the Zariski open subset $Y^{\tr}_\eta$ of $Y_\eta$ where the log
  structure is trivial, and let $\bar \eta_1,\bar \eta_2$ be
  geometric points above them. 
Let $\bar s_2\to\bar s_1$ be a specialization
  of their log reductions, then there is a
  cospecialization homomorphism $\gtempgeom(X_{\bar\eta_1})^{\mbb
    L}\to\gtempgeom(X_{\bar\eta_2})^{\mbb L}$, which is an isomorphism if
  $\bar s_1$ and $\bar s_2$ lie in the same stratum of
  $Y$.\end{thm}
If one does not assume the geometric fibers of $X\to Y$ to have interiorly
free polysimplicial sets, there still is the cospecialization homomorphism
if $\bar s_2$ is the generic point of a stratum.\\

Let us come back to our moduli space of pointed stable elliptic curves with
high enough level
structure $M$ over $O_K$, and let $C$ be the canonical stable elliptic curve on
$M$. If $\eta_1$ and $\eta_2$ are two Berkovich points of $M_\eta$, they
are in $M_\eta^{\tr}$ if and only if $C_{\eta_1}$ and $\C_{\eta_2}$ are
smooth. $C\to M$, endowed with their natural log-structures over
$(O_K,O_K^*)$, is a polystable morphism of log schemes, and $M\to O_K$ is
log smooth. Since the polysimplicial set of a semistable curve over a
separably closed field is interiorly free, the polysimplicial sets of the
geometric fibers of $C\to M$ are interiorly
free, so that one can apply our theorem. One thus get a cospecialization
outer morphism $\gtemp(C_{\overline \eta_1})\to\gtemp(C_{\overline \eta_2})$ for
every specialization $\overline s_2\to\overline s_1$, which is an
isomorphism if $\overline s_1$ and $\overline s_2$ are in the same stratum
of $M_s$. Since the moduli stack
of pointed stable elliptic curves over $\Spec k$ has only two stratum, one
corresponding to smooth elliptic curves $M_0$ and one to singular curves $M_1$,
one gets that $\gtemp(E_1)^{\mbb L}\simeq\gtemp(E_2)^{\mbb L}$ if $E_1$ and $E_2$
are two curves with good reduction or two Tate curves (the isomorphism
depends on choices of cospecializations). Since $M_1$ is in the closure of
$M_0$ one gets a morphism from the tempered fundamental group of a Tate
curve to the tempered fundamental group of an elliptic curve with good
reduction.\\

The first thing we need to construct the cospecialization homomorphism for
tempered fundamental groups is a
specialization morphism for the $(p')$-log geometric fundamental groups of
$X_{\bar s_1}$ and $X_{\bar s_2}$ that extends any $(p')$-log geometric
covering of $X_{s_1}$ to a két neighborhood of $s_1$ (if one has
such a specialization morphism, by comparing it to the fundamental groups
of $X_{\bar\eta_1}$ and $X_{\bar\eta_2}$ and using Grothendieck's
specialization theorem, we will easily get that it must be an
isomorphism). This specialization morphism is easily deduced
from~\cite{org2} if $s_1$ is a strict point of $Y$, \emph{i.e.} the log
structure of $s_1$ is just the pull back of the log structure of $Y$. Thus
we will study the invariance of the log geometric fundamental group by
change of fs base point. Then we have to study the combinatorial behavior
of a két covering with respect to cospecialization. Cospecialization
morphisms of the polysimplicial set of the fibers of a stricly polystable
log fibration are already given by Berkovich in~\cite[cor. 6.2]{berk2} when
$Y$ is plurinodal and $s_1,s_2$ are the topological (\emph{i.e.} not with
value in a field, or equivalently with value in the local field at their
image in $Y$ and with strict log structure) generic points of strata of
$Y$, without any assumption of properness. The construction easily
extends étale locally to our situation if one only still assumes $s_2$ to
be the generic point of a stratum. To get a cospecialization morphism of
geometric polysimplicial sets, we will have to prove that after further két
localization at $\bar s_1$, the strata of $X_{s_2}$ whose closures
meet $X_k$ are geometrically connected. This will follow from the fact,
that after some localization, the closure of those strata are flat over
their image in $Y$ and have reduced geometric fibers. One then descends these
cospecialization morphisms we had étale locally. In the initial proper
case, this cospecialization morphism is an isomorphism if $s_1$ and $s_2$
lie in the same stratum of $Y$ and the polysimplicial set of $X_{\bar s_2}$
is interiorly free. These cospecialization morphisms commute with két
coverings, and thus will give us the wanted cospecialization morphisms.\\

\subsection*{}
Let us now discuss the organization of the paper.\\
The first section of this paper will be devoted to recall the main tools we
will need later. We will recall the definition of the
tempered fundamental group and its basic properties. We will also consider
an $\mbb L$-version of the tempered 
fundamental group, where $\mbb L$ is a set of prime numbers ($\mbb
L$-tempered fundamental groups were already introduced in~\cite{mochi} in
the case of curves). We will then recall the basics of log geometry, especially
the theory of két coverings and log fundamental groups. We will end this
part by recalling the topological structure of the Berkovich space of a
pluristable formal scheme, as studied in~\cite{berk2} and in~\cite{berk3}.\\

In §2, we define the tempered fundamental group of a pluristable
log scheme $X$ over a log point. To do this, we define a functor $\C$ from
the Kummer étale site of our pluristable log scheme $X$ to the category of
polysimplicial sets (which extends the definition of the polysimplicial set
associated to a pluristable scheme defined by Berkovich
in~\cite{berk2}). Thus, for any Galois két covering $Y$ of $X$, there is an
action of $\Gal(X/Y)$ on $\C(Y)$, which defines an 
extension of $\Gal(X/Y)$ by $\gtop(|\C(Y)|)$. The tempered fundamental
group of $X$ will then be defined to be the projective limits of these
extensions, when $Y$ runs through pointed két Galois coverings of $X$. As
for the tempered fundamental group of a Berkovich space, one also defines
$\mbb L$-versions of 
the tempered fundamental group of our pluristable log scheme. One also defines
a log geometric version by taking the projective
limit under connected két extensions of the log point.\\

In §3, for a proper, generically smooth and pluristable scheme $X$ over a
complete discretely valued ring $O_K$ (thus endowed with a canonical log
structure), we construct a specialization morphism between the $\mbb 
L$-tempered fundamental group of the generic fiber, considered as a
Berkovich space, and the $\mbb L$-tempered fundamental group of the special
fiber endowed with the inverse image log structure, which is an isomorphism
 if the residual characteristic of $K$ is not in $\mbb L$.\\ 
This specialization morphism is induced by the specialization morphism from
the algebraic fundamental group of the generic fiber to the log fundamental
group of the special fiber, and by the fact that the geometric realization
of the polysimplicial set $|\C(Y)|$ of a két covering of the special fiber
of $X$ is canonically homotopically equivalent to the Berkovich space
$Y^{\an}_\eta$ of the corresponding étale covering of the generic
fiber. This homotopy equivalence is obtained by extending the strong
deformation retraction of $X^{\an}_\eta$ to a strong deformation retraction
of $Y^{\an}_\eta$ onto a subset canonically homeomorphic to $|\C(Y)|$.\\

In §4, we start by studying specialization of the $(p')$-log fundamental group of
a proper, log smooth and saturated morphism of fs log schemes.\\
We then construct cospecialization morphisms between the
polysimplicial sets of the fibers of a strictly polystable fibration over a log
regular Zariski log scheme. This cospecialization question is already studied
in~\cite{berk2}, when the base scheme is strictly plurinodal for the
generic points of the strata. Our cospecialization morphisms extend to
cospecialization morphisms of the geometric polysimplicial sets of the
fibers of a két covering of a strictly polystable fibration.\\
The two cospecialization theories fit together, and thus we obtain
cospecialization morphisms between the $(p')$-geometric tempered
fundamental groups of the fibers of our strictly polystable log fibration.\\
Thanks to the isomorphisms between the $(p')$-geometric tempered
fundamental group of the fiber over a discretely valued Berkovich point of the generic part of our base log
scheme and the $(p')$-geometric tempered fundamental group of the fiber
over the reduction log point, we will get theorem~\ref{sp}.\\

This work is part of a PhD thesis. I would like to thank my advisor, Yves André, for
suggesting me to work on the cospecialization of the tempered fundamental
group and taking the time of reading and correcting this work. I would also
like to thank Luc Illusie and Fumiharu Kato for taking interest in my
problem about the invariance of geometric log fundamental groups by base change.\\

\section{Preliminaries}

\subsection{Tempered fundamental group}

Let $K$ be a complete nonarchimedean field.\\
Let $\mbb L$ be a set of prime numbers (for example, we will denote by
$(p')$ the set of all primes except the residual characteristic $p$ of
$K$). An $\mbb L$-integer will be an integer which is a product of elements
of $\mbb L$.\\
Following~\cite[§4]{andre2}, a $K$-\emph{manifold} will be a smooth paracompact strictly $K$-analytic space in the sense of Berkovich. For example, if $X$ is a smooth algebraic $K$-variety, $X^{\an}$ is a $K$-manifold (and in fact, we will mainly be interested in those spaces). Then, thanks to~\cite{berk2}, any $K$-manifold is locally contractible (we will explain in more detail the results of~\cite{berk2} in section~\ref{berkspaces}). In particular, it has a universal covering.\\
A morphism $f:S'\to S$ is said to be an \emph{étale covering} if $S$ is covered by open subsets $U$ such that $f^{-1}(U)=\coprod V_j$ and $V_j\to U$ is étale finite~(\cite{dJ1}).\\
For example, étale $\mbb L$-finite coverings (\emph{i.e.} étale finite
coverings such that the order of every connected component is an $\mbb L$-integer), also called \emph{$\mbb
  L$-algebraic coverings}, and coverings in the usual topological sense for
the Berkovich topology, also called \emph{topological coverings}, are
étale coverings.\\
Then, André defines tempered coverings in \cite[def. 2.1.1]{andre1}. We generalize this definition to
$\mbb L$-tempered coverings as follows:
\begin{dfn} \label{def:rvt:temp}
An étale covering $S' \to S$ is \emph{$\mbb L$-tempered} if it is a
 quotient of the composition of a topological covering $T' \to T$ and of a
 $\mbb L$-finite étale covering
 $T \to S$.
\end{dfn}
This is equivalent to say that it becomes a topological covering after
pullback by some $\mbb L$-finite étale covering.\\
We denote by $\Covtemp(X)^{\mbb L}$ (resp. $\Covalg(X)^{\mbb L}$,
$\Covtop(X)$) the category of $\mbb L$-tempered coverings (resp. $\mbb L$-algebraic coverings, topological coverings) of $X$ (with the obvious morphisms).\\

A geometric point of a $K$-manifold $X$ is a morphism of Berkovich spaces $\mcal M(\Omega)\to X$ where $\Omega$ is an algebraically closed complete isometric extension of $K$.\\
Let $\bar x$ be a geometric point of $X$. Then one has a functor \[F^{\mbb L}_{\bar
x}:\Covtemp(X)^{\mbb L}\to\Set\] which maps a covering $S\to X$ to the set
$S_{\bar x}$. If $\bar x$ and $\bar x'$ are two geometric points, then
$F^{\mbb L}_{\bar x}$ and $F^{\mbb L}_{\bar x'}$ are (non canonically) isomorphic (\cite[prop. 2.9]{dJ1}).\\
The tempered fundamental group of $X$ pointed at $\bar x$ is
\[\gtemp(X,\bar x)^{\mbb L}=\Aut F^{\mbb L}_{\bar x}.\]
When $X$ is a smooth algebraic $K$-variety, $\Covtemp(X^{\an})^{\mbb L}$
and $\gtemp(X^{\an},\bar x)^{\mbb L}$ will also be denoted simply by
$\Covtemp(X)^{\mbb L}$ and $\gtemp(X,\bar x)^{\mbb L}$.\\
By considering the stabilizers $(\Stab_{F(S)}(s))_{S\in \Covtemp(X)^{\mbb
L},s\in F_{\bar x}(S)}$ as a basis of open subgroups of $\gtemp(X,\bar
x)^{\mbb L}$, $\gtemp(X,\bar x)^{\mbb L}$ becomes a topological group. It is a prodiscrete topological group.\\
When $X$ is algebraic, $K$ of characteristic zero and has only countably
many finite extensions in a fixed algebraic closure $\overline K$,
$\gtemp(X,\bar x)^{\mbb L}$ has a countable fundamental system of neighborhood of $1$ and all its discrete quotient groups are finitely generated~(\cite[prop. 2.1.7]{andre1}).\\

If $\bar x$ and $\bar x'$ are two geometric points, then $F^{\mbb L}_{\bar
x}$ and $F^{\mbb L}_{\bar x'}$ are (non canonically) isomorphic
(\cite[prop. 2.9]{dJ1}). Thus, as usual, the tempered fundamental group
depends on the basepoint only up to inner automorphism. This topological
group, considered up to conjugation, will sometimes be denoted simply
$\gtemp(X)^{\mbb L}$.\\
The full subcategory of tempered coverings $S$ for which $F^{\mbb L}_{\bar
x}(S)$ is $\mbb L$-finite is equivalent to $\Covalg(S)^{\mbb L}$, hence
\[\big(\gtemp(X,\bar x)^{\mbb L}\big)^{\mbb L}=\ga(X,\bar x)^{\mbb L}\] (where $(\ )^{\mbb L}$ denotes the pro-$\mbb L$ completion).\\
For any morphism $X\to Y$, the pullback defines a functor
$\Covtemp(Y)^{\mbb L}\to\Covtemp(X)^{\mbb L}$. If $\bar x$ is a geometric
point of $X$ with image $\bar y$ in $Y$, this gives rise to a continuous
homomorphism \[\gtemp(X,\bar x)^{\mbb L}\to\gtemp(Y,\bar y)^{\mbb L}\]
(hence an outer morphism $\gtemp(X)^{\mbb L}\to\gtemp(Y)^{\mbb L}$).\\
One has the analog of the usual Galois correspondence:
\begin{thm}[{\cite[th. 1.4.5]{andre1}}]  \label{galcorr} $F^{\mbb
L}_{\bar x}$ induces an equivalence of categories between the category of
direct sums of $\mbb L$-tempered coverings of $X$ and the category
$\gtemp(X,\bar x)^{\mbb L}\tSet$ of discrete sets endowed with a
continuous left action of $\gtemp(X,\bar x)^{\mbb L}$.\end{thm}

If $S$ is a $\mbb L$-finite Galois covering of $X$, its universal
topological covering $S^{\infty}$ is still Galois and every connected $\mbb
L$-tempered covering is dominated by such a Galois $\mbb L$-tempered covering.\\
If $((S_i,\bar s_i))_{i\in \mbf N}$ is a cofinal projective system (with
morphisms $f_{ij}:S_i\to S_j$ which maps $s_i$ to $s_j$ for $i\geq j$) of
geometrically pointed Galois $\mbb L$-finite étale coverings of $(X,\bar
x)$, let $((S^{\infty}_i,\bar s^{\infty}_i))_{i\in\mbf N}$ be the
projective system, with morphisms $f_{ij}^\infty$ for $i\geq j$, of its
pointed universal topological coverings. Then $F^{\mbb L}_{\bar
x}(S^{\infty}_i)=\gtemp(X,\bar x)^{\mbb L}/\Stab_{F(S^\infty_i)}(\bar
s^\infty_i)$ is naturally a quotient group $G$ of $\gtemp(X,\bar x)^{\mbb
L}$  for which $s^{\infty}_i$ is the neutral element. Moreover $G$ acts by
$G$-automorphisms on $F^{\mbb L}_{\bar x}(S^{\infty}_i)$ by right
translation (and thus on $S^{\infty}_i$ thanks to the Galois
correspondence~(theorem \ref{galcorr})). Thus one gets a morphism
$\gtemp(X,\bar x)^{\mbb L}\to\Gal(S^{\infty}_i/X)$. As $f_{ij}^\infty(s_i^\infty)=s_j^\infty$, these morphisms are compatible with $\Gal(S^{\infty}_i/X)\to \Gal(S^{\infty}_j/X)$.\\
Then, thanks to~\cite[lem. III.2.1.5]{andre1},
\begin{prop}\label{limproj} \[\gtemp(X,\bar x)^{\mbb L}\to\varprojlim \Gal(S^{\infty}_i/X)\] is an isomorphism.\end{prop}

\subsection{Log fundamental groups}\label{loggeometry}
This part is a reminder of the theory of log schemes, as can be found
in~\cite{kato} and \cite{ogus}, and of the log
fundamental groups, as can be found in~\cite{ill} or \cite{stix}.
\subsubsection{Log schemes}
All monoids here are commutative with units, and morphisms of monoids map
the unit to the unit. If $P$ is a monoid, $P^{\gp}$ will be its group
envelope, $P^*$ the group of invertible elements of $P$ and $\overline
P=P/P^*$.\\
If $a,b\in P$, we will write $a|b$ if there is $c\in P$ such that $b=ac$.\\
A monoid $P$ is \emph{sharp} if $P^*=\{1\}$.\\ 
A monoid $P$ is \emph{integral} if the morphism $P\to P^{\gp}$ is
injective. A monoid is \emph{fine} if it is integral and finitely generated.
An integral monoid is \emph{saturated} if $a\in P^{\gp}$ is in $P$
if there exists $n$ such that $a^n\in P$.\\
If $P$ is a fine and saturated monoid (or fs for short),
$P^{\gp}/P^*=\overline P^{\gp}$ is a
free abelian group of finite type, and thus there is a (non canonical) section $\overline
P^{\gp}\to P^{\gp}$ and it induces a decomposition $P=P^*\oplus \overline P$.\\
A morphism $f:P\to Q$ of monoids is \emph{local} if $f^{-1}(Q^*)=P^*$
A morphism $P\to Q$ of integral monoids is \emph{exact} if $P$ is the inverse
image of $Q$ in $P^{\gp}$.\\

A \emph{prime} $\fk p$ of an fs monoid $P$ is a subset of $P$ such that if
$p\in\fk p$ and $p'\in P$ then $p+p'\in\fk P$, and if $p,p'\in P$ and
$p+p'\in\fk p$ then $p\in\fk p$ or $p'\in \fk p$.\\
A subset $F$ of an fs monoid $P$ is called a \emph{face} if $P\backslash F$ is a
prime (in particular $F$ is a submonoid of $P$).\\
$\Spec P$ denotes the topological space of primes of $P$, where $(D(f)=\{\fk
p,f\notin\fk p\})_{f\in P}$ is a basis of the topology of $\Spec P$.\\
If $f:P\to Q$ is a morphism of fs monoids and $\fk q$ is a prime of $Q$,
then $f^{-1}(\fk q)$ is a prime of $P$, hence a map $\Spec Q\to\Spec P$.\\

A \emph{pre-log structure} on a scheme $X$ is a pair $(M,\alpha)$ where $M$ is a
sheaf of monoids on $X_{\et}$, and $\alpha:M\to (O_X,.)$ is a morphism of
sheaves of monoids, where $O_X$ is the canonical sheaf of $X$ and $.$
is the multiplication on $O_X$. A pre-log structure will be a \emph{log structure} if the induced map
$\alpha^{-1}(O_X^*)\to O_X^*$ is an isomorphism. A \emph{log scheme} $X$ is a
scheme (the \emph{underlying scheme} $\mring{X}$ of the log scheme) with a log structure on it.\\
The forgetful functor from log structures on $X$ to pre-log structures on $X$ admits a
left adjoint $(M,\alpha)\mapsto (M^a,\alpha^a)$ where $M^a$ is the
amalgamated sum of $Q$ and $O_X$ along $\alpha^{-1}(O^*_X)$ (this log
structure is called the log structure \emph{associated} to $(M,\alpha)$).\\
A \emph{morphism of log schemes} $f:(X,M,\alpha)\to (Y,N,\beta)$ is a morphism of schemes
$f:X\to Y$ with a morphism of sheaves of monoids $f^{-1}N\to M$ compatible
with $\alpha$ and $\beta$. Then $f^{-1}N\to M$ is necessarily a local
morphism of sheaves of monoids.\\
A log scheme $X$ is \emph{integral} if for every geometric point $\bar x$ of
$\mring{X}$, $M_{\bar x}$ is integral.\\
If $Y=(Y,N,\beta)$ is a log scheme and $X\to\mring{Y}$ is a morphism of schemes, the log
structure on $X$ associated to $(f^{-1}N,f^{-1}\beta)$ is called the
\emph{inverse image log structure} and is denoted $f^*N$. A morphism of log
schemes $f:(X,M,\alpha)\to (Y,N,\beta)$ is \emph{strict} if the
induced morphism $f^*N\to M$ is an isomorphism (if $X$ is integral,
this is equivalent to say that
$\overline N_{f(\bar x)}\to\overline M_{\bar x}$ is an isomorphism for
every geometric point $\bar x$ of $X$).\\

If $P$ is a monoid the log structure on $\Spec \mbf Z[P]$ associated to the
pre-log structure defined by $P\to\mbf Z[P]$ is called the canonical log
structure. There is a canonical morphism $\Spec \mbf Z[P]\to\Spec P$ which
maps a prime ideal $I$ of $\mbf Z[P]$ to $I\cap P$.\\
A (global) \emph{chart} modeled on a monoid $P$ of a log scheme $X$ is a morphism
from the constant sheaf $P_X\to M_X$ inducing an isomorphism on the
associated log structures. This also amounts to giving a strict morphism
$X\to\Spec \mbf Z[P]$, where $\mbf Z[P]$ is endowed with its canonical log
structure.\\
If $\bar x$ is a geometric point, an integral chart $P\to M$ is \emph{exact at $\bar
  x$} (resp. \emph{good at $\bar x$}) if $\overline P\to \overline
M_{\bar x}$ (resp. $P\to\overline M_{\bar x}$) is an isomorphism.\\
A log scheme is \emph{fine} (resp. \emph{fine and saturated} or fs for
short) if it is integral and, locally for the étale topology, it admits a chart modeled on a
finitely generated and integral (resp. finitely generated and saturated)
monoid.\\
We will mainly work in the category of fs log schemes. There are fiber
products in this category, but in general, taking the underlying scheme
does not commute with fiber products.\\
If $X\to\Spec P$ is an fs chart and $\overline x$ is a geometric point of
$X$ which maps to $\fk p\in\Spec P$ and let $F=P\backslash \fk p$, then
$P\to M_{\bar x}$ induces an isomorphism $\overline{F^{-1}P}\to\overline
M_{\bar x}$. Moreover $\Spec \mbf Z[F^{-1}P]\to\Spec\mbf Z[P]$ is an open
embedding corresponding to the preimage of $D(\fk p)=\{\fk p'|\fk p\subset
\fk p'\}\subset \Spec P$. Thus $F^{-1}P$ induces an exact chart of a
Zariski neighborhood of $\overline x$. But one can then choose some
decomposition $F^{-1}P=\overline{F^{-1}P}\oplus(F^{-1}P)^*$, and the
induced morphism $\overline{F^{-1}P}\to F^{-1}P\to M_X$ is a good chart at
$\overline x$.\\

Sometimes, we may have to use log structures on the Zariski site. Let
$\epsilon:X_{\Zar}\to X_{\et}$ the natural projection. We will say that
a log structure $M$ on $X$ is \emph{Zariski} (and the log scheme $X$ is
\emph{log Zariski}) if $\epsilon^*\epsilon_*M\to M$ is an
isomorphism. If $X$ is an fs log scheme, the log structure is Zariski if
and only if it has fs charts locally on the Zariski topology. In particular
any fs log scheme is étale locally log Zariski.\\

If $f:X\to Y$ is a morphism of fine log schemes, a \emph{chart} of $f$ is
given by a chart $X\to\Spec \mbf Z[P]$, a chart $Y\to\Spec \mbf Z[Q]$ and a
morphism $Q\to P$ such that the corresponding square of log schemes
commute. Any morphism of fine log schemes has charts étale locally.

A morphism of fine log schemes $f:X\to Y$ is \emph{log smooth}
(resp. \emph{log étale}) if étale locally on $X$ and $Y$, $f$ admits a
chart $u:Q\to P$ such that the kernel and the torsion part of the cokernel
(resp. the kernel and the cokernel) of $u^{\gp}$ are finite groups of order
invertible on $X$ and $X\to Y\times_{\Spec \mbf Z[Q]}\Spec \mbf Z[P]$ is
étale.\\
There are valuative characterizations of log étale and log smooth
morphisms. Log étale and log smooth morphisms are stable under base change
and composition.

A morphism $h:Q\to P$ of fs monoids is \emph{Kummer} (resp. $\mbb
L$-Kummer)
if $h$ is injective
and for every $a\in P$, there exists an integer (an $\mbb L$-integer) $n$ such that $na\in h(Q)$
(note that if $Q\to P$ is Kummer, $\Spec P\to\Spec Q$ is an homeomorphism).\\
A morphism $f:X\to Y$ of fs log schemes is said to be \emph{Kummer} (resp. \emph{exact})
if for
every geometric point $\bar x$ of $X$, $\overline M_{Y,f(\bar x)}\to \overline
M_{X,\bar x}$ is Kummer (resp. exact).\\
A Kummer morphism $X\to Y$ of fs log schemes is a \emph{Kummer universal
  homeomorphism} (or \emph{kuh} for short)
if the underlying map of schemes is an homeomorphism after any fs base
change (\cite[def. 2.1]{vidal}).\\
A Kummer morphism $q:X\to Y$ is kuh if and only if $\mring q$ is a universal
homeomorphism (\emph{i.e.} is integral, radicial and surjective) and for
any geometric point $\bar x$ of $X$, $\overline M_{Y,q(\bar x)}=\overline
M_{X,\bar x}$ is $p$-Kummer, where $p$ is the residual characteristic of
$\bar x$ (\cite[thm. 2.7]{vidal}).\\
For example, if $P\to Q$ is a $p$-Kummer morphism of fs monoid, and $A$ is a ring of characteristic $p$,
then $\Spec A[Q]\to\Spec A[P]$ is kuh.\\

An fs log scheme $X$ is \emph{log regular} if for every geometric point
$\bar x$, $O_{X,\bar x}/I_{\bar x}O_{X,\bar x}$ is regular and
$\dim(O_{X,\bar x})=\dim(O_{X,\bar x}/I_{\bar x}O_{X,\bar x})+\rk(\overline
M_{\bar x}^{\gp})$ where $I_{\bar x}$ is the ideal of $O_{X,\bar x}$
generated by the image of $M_{\bar x}\backslash O_{X,\bar x}^*$
(\cite[def. 2.2]{niziol}).\\
Log Zariski log regular log schemes are studied in~\cite{kato2}.\\
If $X$ is log regular, $\mring X$ is normal. The subset
$X_{\tr}=\{x,\overline M_{\bar x}=\{1\}\}$ of $X$ is a dense open
subset of $X$ and \[M=O_X\cap j_*O_{X_{\tr}}^*\] where $j$ is the open
embedding $X_{\tr}\to X$~(\cite[prop. 2.6]{niziol}).\\
If $Y$ is log regular and $X\to Y$ is log smooth, then $X$ is log
regular (\cite[thm. 8.2]{kato2}).

\subsubsection{Két coverings}

A morphism of fs log scheme is \emph{Kummer étale} (or két for short) if it
is Kummer and log étale.\\
A morphism $f$ is két if and only if étale locally it is deduced by strict
base change and étale localization from a map $\Spec \mbf Z[P]\to\Spec \mbf
Z[Q]$ induced by a Kummer map $Q\to P$ such that $nP\subset Q$ for some $n$
invertible on $X$.\\
In fact if $f:Y\to X$ is két, $\bar y$ is a geometric point of $Y$, and
$P\to M_X$ is an exact chart of $X$ at $f(\bar y)$, there is an étale
neighborhood $U$ of $\bar x$ and a Zariski open neighborhood $V\subset
f^{-1}(U)$ of $\bar y$ such that $V\to U$ is isomorphic to $U\times_{\Spec
  \mbf Z[P]}\Spec\mbf Z[Q]$ with $P\to Q$ a $\mbb L$-Kummer morphism where
$\mbb L$ is the set of primes invertible on $U$ (\cite[Prop. 3.1.4]{stix}).\\
Két morphisms are open and quasi-finite.\\

The category of két fs log schemes over $X$ (any $X$-morphism
between two such fs log schemes is then két) where the covering families
$(T_i\to T)$ of $T$ are the families that are set-theoretical covering
families (being a set-theoretical covering két family is stable under fs base change) is a site. We
will denote by $X_{\ket}$ the corresponding topos.\\
Any locally constant finite object of $X_{\ket}$ is representable. A két fs
log scheme over $X$ that represents such a locally constant finite sheaf
will be called a \emph{két covering} of $X$. We will
denote by $\KCov(X)$ the category of két coverings
of $X_{\ket}$.\\
A log geometric point is a log scheme $s$ such that $\mring{s}$ is the
spectrum of a
separably closed field $k$ such that $M_s$ is saturated and multiplication
by $n$ on $\overline
M_s$ is an isomorphism for every $n$ prime to the characteristic of $k$.\\
A log geometric point of $X$ is a morphism $x:s\to X$ of log schemes where $s$ is a log
geometric point. A két neighborhood $U$ of $s$ in $X$ is a morphism $s\to
U$ of $X$-log schemes where $U\to X$ is két. Then if $x$ is a log geometric
point of $X$, the functor $F_x$ from $X_{\ket}$ to $\Ens$ defined by $\mcal F\mapsto
\varinjlim_U \mcal F(U)$ where $U$ runs through the directed category of két
neighborhoods of $x$ is a point of the topos $X_{\ket}$ and any point of
this topos is isomorphic to $F_x$ for some log geometric point and they
make a conservative system.\\
 One also defines the log strict
  localization $X(x)$ to be the inverse limit in the category of saturated
log schemes of the két neighborhoods of $x$. If $x$ and $y$ are log
geometric points of $x$, a \emph{specialization} of log geometric points $x\to y$
is a morphism $X(x)\to X(y)$ over $X$.\\
A specialization $x\to y$ induces a canonical morphism  $F_x\to F_y$ of functors.\\
 If there is a specialization $x\to
y$ of the underlying topological points, then there is some specialization $x\to y$ of log
geometric points.\\
If $X$ is connected, for any log geometric point $x$ of $X$, $F_x$ induces
a fundamental functor $\KCov(X)\to\fSet$ of the Galois category
$\KCov(X)$.\\
One then denotes by $\glog(X,x)$ the profinite group of automorphisms of
this fundamental functor.\\

Strict étale surjective morphisms satisfy effective descent for két
coverings (\cite[prop. 3.2.19]{stix}).\\
If $f:S'\to S$ is an exact morphism of fs log schemes such that $\mring{f}$
is proper, surjective and of finite presentation, then $f$ satisfies
effective descent for két coverings (\cite[th. 3.2.25]{stix}).\\

If $X$ is a log regular fs log scheme, and all the primes of $\mbb L$ are
invertible on $X$, then $\KCov(X)^{\mbb L}\to\Covalg(X_{\tr})^{\mbb L}$ is
an equivalence of categories (\cite[th. 7.6]{ill}).\\
\begin{thm}[{\cite[cor. 2.3]{org2}}]\label{orgsp} Let $S$ be a strictly local scheme
  with closed point $s$, and let $X$ be a connected fs log scheme such that
  $\mring X$ is proper over $S$. Then
\[\KCov(X)\to\KCov(X_s)\]is an equivalence of categories.\end{thm}
If $q:X\to Y$ is kuh, $q^*:Y_{\ket}\to X_{\ket}$ is an equivalence of
categories (\cite[th. 0.1]{vidal})

\subsubsection{Saturated morphisms}

A morphism of fs monoids $P\to Q$ is \emph{integral} if, for any morphism
of integral monoids $P\to Q'$, the amalgamated sum $Q\oplus_PQ'$ is still
integral.\\
A integral morphism of fs monoids $P\to Q$ is \emph{saturated} if, for any morphism
of fs monoids $P\to Q'$, the amalgamated sum $Q\oplus_PQ'$ is still
fs.\\
If $f:P\to Q$ is saturated and $F'$ is a face of $Q$,
$\overline{F^{-1}P}\to\overline{{F'}^{-1}Q}$ is also saturated where $F=f^{-1}(F')$.\\
\begin{lem}\label{satface}If $\phi:P\to Q$ is an integral (resp. saturated) morphism of fs monoids and $F'$ is a face of $Q$, let
$F=\phi^{-1}(F')$. Then $F\to F'$ is also integral (resp. saturated).\end{lem}
\dem
To prove that $F\to F'$ is integral, thanks to~\cite[prop. I.4.3.11]{ogus},
one only has to prove that if $f'_1,f'_2\in F'$ and $f_1,f_2\in F$ are such
that $f'_1\phi(f_1)=f'_2\phi(f_2)$, there are $g'\in F'$ and $g_1,g_2\in
F$ such that $f'_1=g'\phi(g_1)$ and $f'_2=g'\phi(g_2)$.\\
But there exists $g'\in Q$ and $g_1,g_2\in P$ that satisfies those
properties since $P\to Q$ is integral. But, since $F'$ is a face of $Q$,
$g',\phi(g_1),\phi(g_2)$ must be in $F'$, and thus $g_1$ and $g_2$ are in
$F$.\\
Thanks to a criterion of T. Tsuji, an integral morphism of fs monoids $f:P_0\to Q_0$ is saturated if and only
if for any $a\in P_0, b\in Q_0$ and any prime number $p$ such that
$f(a)|b^p$, there exists $c\in P_0$ such that $a|c^p$ and $f(c)|b$. Let
$a\in F, b\in F'$ and $p$ be a prime such that $\phi(a)|b^p$. Then since
$\phi:P\to Q$ is saturated, there exists $c\in P$ such that $a|c^p$ and
$f(c)|b$. But $f(c)|b$ implies that $f(c)\in F'$, whence $c\in F$.
\findem

A morphism $f:Y\to X$ of fs log schemes is \emph{saturated} if for any geometric
point $\bar y$ of $Y$, $\bar M_{X,f(\bar y)}\to\bar M_{Y,\bar y}$ is
saturated.\\
If $Y\to X$ is saturated and $Z\to X$ is a morphism of fs log schemes, then
the underlying scheme of $Z\times_XY$ is $\mring Z\times_{\mring X}\mring
Y$.\\
 
If $P\to Q$ is a local and integral (resp. saturated) morphism of fs
monoids and $P$ is sharp, the morphism $\Spec \mbf
Z[Q]\to\Spec\mbf Z[P]$ is flat (resp. separable, \emph{i.e.} flat
with geometrically
reduced fibers, cf.~\cite[cor. 4.3.16]{ogus} and \cite[rem. 6.3.3]{satmorph}).\\
Let $f:X\to Y$ be log smooth and $\bar x$ be a geometric point of
$X$. Etale locally on $Y$, there is a good chart 
$Y\to\Spec P$ at $\bar y$. Then, thanks to~\cite[prop. A.3.1.1]{satmorph},
there is étale locally at $x$ a chart $P\to Q$ of $Y\to X$ such that $Y\to
\Spec \mbf Z[Q]\times_{\mbf Z[P]}X$ is étale and $X\to\Spec Q$ is exact at
$x$. Thus if $f$ is integral (resp. saturated), $P\to Q$ is a local and
integral (resp. saturated) morphism of fs monoids and $P$ is sharp. Thus
$f$ is flat (resp. separable).\\
 
If $P\to Q$ is an integral morphism of fs monoids, there exists an integer
$n$ such that the pullback $P_n\to Q'$ of $P\to Q$ along
$P\stackrel{n}{\to}P=P_n$ is saturated (theorem~\cite[A.4.2]{satmorph}).\\
Moreover if $P\to Q$ factors through $Q_0$ such that $P\to Q_0$ is saturated
and $Q_0\to Q$ is $\mbb L$-Kummer, $n$ can be chosen to be an $\mbb L$-integer.\\
Thus, if $Z''\to Z'$ is a két covering and $Z'\to Z$ is saturated, log smooth
and proper, then for any
log geometric point $z$ of $Z$ there is a két neighborhood $U$ of $z$ such
that $Z''_U\to U$ is saturated (and all the properness assumptions can be
removed and replaced by the quasicompactness of $Z''$ if $Z$ is just an fs
log point, \emph{i.e.} its underlying scheme is the spectrum of a field).

\subsubsection{Specialisation of log fundamental groups}
Let $X\to S$ be a proper and saturated morphism of log schemes,
  and let $Y\to X$ be a két covering. Let $s$ and $s'$ be two points of $S$
  and assume that one has a specialization $\bar
    s'\to\bar{s}$ (where $\bar s$ and
  $\bar s'$ are some log geometric points over $s$ and
  $s'$).\\
Let $Z$ be the strictly local scheme of $S$ at $s$ endowed with the inverse
image log structure, and let $z$ be its closed point, endowed with the
inverse image log structure.\\
One has the following arrows (defined up to inner homomorphisms):
\[\ggeom(Y_{s}/s)^{(p')}\to\ggeom(Y_z/z)^{(p')}\stackrel{\simeq}{\to}\ggeom(Y_Z/Z)^{(p')}\leftarrow\ggeom(Y_{s'}/s')^{(p')}\]
where the second homomorphism are isomorphisms according to
theorem~\ref{orgsp}. The first one is an isomorphism according
to~\cite[cor. 2.6]{cosp1}.\\
\begin{thm}[{\cite[cor. 2.7]{cosp1}}]\label{logsp}One has a specialization morphism
\[\ggeom(Y_{s'}/s')^{(p')}\to\ggeom(Y_{s}/s)^{(p')}\] that factors through $\ggeom(Y_Z/Z)^{(p')}$.\end{thm}

\subsection{Skeleton of a Berkovich space with pluristable reduction}\label{berkspaces}

Let $K$ be a complete nonarchimedean field and let $O_{K}$ be its ring of integers.\\
If $\fk X$ is a locally finitely presented formal scheme over $O_K$, $\fk
X_\eta$ will denote the generic fiber of $\fk X$ in the sense of Berkovich
(\cite[section 1]{berkvc1}).\\ 
Recall the definition of a  polystable morphism of formal schemes:
\begin{dfn}[{\cite[def. 1.2]{berk2}, \cite[section 4.1]{berk3}}] Let
  $\phi:\fk Y\to \fk X$ be a locally finitely presented morphism of formal
  schemes over $O_K$. 
\begin{enumerate}[(i)]
\item $\phi$ is said to be \emph{strictly polystable} if, for every point
  $y\in\fk Y$, there exists an open affine neighborhood $\fk X'=\Spf(A)$ of
  $x:=\phi(y)$ and an open neighborhood $\fk Y'\subset\phi^{-1}(\fk X')$ of
  $x$ such that the induced morphism $\fk Y'\to \fk X'$ factors through an
  étale morphism $\fk Y'\to\Spf(B_0)\times_{\fk X'}\cdots\times_{\fk
    X'}\Spf(B_p)$ where each $B_i$ is of the form
  $A\{T_0,\cdots,T_{n_i}\}/(T_0\cdots T_{n_i}-{a_i})$ with $a\in A$ and $n\geq
  0$. It is said to be \emph{nondegenerate} if one can choose $X'$, $Y'$
  and $(B_i,a_i)$ such that $\{x\in(\Spf(A)_\eta)|a_i(x)=0\}$ is nowhere
  dense. 
\item $\phi$ is said to be \emph{polystable} if there exists a surjective
  étale morphism $\fk Y'\to\fk Y$ such that $\fk Y'\to \fk X$ is strictly
  polystable. It is said to be \emph{nondegenerate} if one can choose $\fk
  Y'$ such that $\fk Y'\to \fk X$ is
  nondegenerate. \end{enumerate}\end{dfn} 
Then a (\emph{nondegenerate}) \emph{polystable fibration} of length $l$
over $\fk S$ is a sequence of (nondegenerate) polystable morphisms
$\underline{\fk X}=(\fk X_l\to\cdots\to\fk X_1\to\fk S)$.\\ 
Then $K\text{-}\mcal Pstf_l^{\et}$ will denote the category of polystable
fibrations of length $l$ over $O_K$, where a morphism $\underline{\fk X}\to
\underline{\fk Y}$ is a collection of étale morphisms $(\fk X_i\to\fk
Y_i)_{1\leq i\leq l}$ which satisfies the natural commutation
assumptions.\\ 
$Pstf_l^{\et}$ will denote the category of couples $(\underline{\fk
  X},K_1)$ where $K_1$ is a complete non archimedean field and
$\underline{\fk X}$ is a polystable fibration over $O_{K_1}$, and a
morphism $(\underline{\fk X},K_1)\to (\underline{\fk Y},K_2)$ is a couple
$(\phi,\psi)$ where $\phi$ is an isometric extension $K_2\to K_1$ and
$\psi$ is a morphism $\underline{\fk X}\to\underline{\fk
  Y}\otimes_{O_{K_2}}O_{K_1}$ in $K_1\text{-}\mcal Pstf_l^{\et}$.\\  

Let $k$ be a field.\\
Let $X$ be a $k$-scheme locally of finite type.\\
The normal locus $\Norm(X^{\red})$ is a dense open subset of $X$. Let us define inductively
$X^{(0)}=X^{\red}$, $X^{(i+1)}=X^{(i)}\backslash \Norm(X^{(i)})$. The
irreducible components of $X^{(i)}\backslash X^{(i+1)}$ are called the
strata of
$X$ (of rank $i$). This gives a partition of $X$. The set of the generic
points of the strata of $X$ is denoted by $\Str(X)$ (There is a natural
bijection with the set of strata of $X$).\\
Berkovich defines another filtration $X=X_{(0)}\subset X_{(1)}\subset\cdots$
such that $X_{(i+1)}$ is the closed subset of points contained in at least
two irreducible components of $X_{(i)}$. $X$ is said to be
\emph{quasinormal} if all of the irreducible components of each $X_{(i)}$,
endowed with the reduced subscheme structure, is normal (this property is
local for the Zariski topology and remains true after étale morphisms). If
$X$ is quasinormal, then $X_{(i)}=X^{(i)}$. $X$ is quasinormal if and only
if the closure of every stratum is normal.\\  
There is a natural partial order on $\Str(X)$  defined by $x\leqslant y$ if
and only if $y\in
\overline{\{x\}}$.\\

Berkovich defines \emph{polysimplicial sets} in~\cite[section 3]{berk2} as
follows.\\ 
For an integer $n$, let $[n]$ denote $\{0,1,\cdots,n\}$.\\ 
For a tuple $\mbf n=(n_0,\cdots,n_p)$ with either $p=n_0=0$ or $n_i\geq 1$
for all $i$, let $[\mbf n]$ denote the set $[n_0]\times\cdots\times [n_p]$
and $w(\mbf n)$ denote the number $p$.\\ 
Berkovich defines a category $\mbf \Lambda$ whose objects are $[\mbf n]$
and morphisms are maps $[\mbf m]\to [\mbf n]$ associated with triples
$(J,f,\alpha)$, where: 
\begin{itemize}
\item $J$ is a subset of $[w(\mbf m)]$ assumed to be empty if $[\mbf m]=[0]$, 
\item $f$ is an injective map $J\to [w(\mbf n)]$,
\item $\alpha$ is a collection $\{\alpha_l\}_{0\leq l\leq p}$, where
  $\alpha_l$ is an injective map $[m_{f^{-1}(l)}]\to [n_l]$  if
  $l\in\Image(f)$, and $\alpha_l$ is a map $[0]\to [n_l]$
  otherwise.\end{itemize} 
The map $\gamma:[\mbf m]\to [\mbf n]$ associated with $(J,f,\alpha)$ takes
$\mbf j=(j_0,\cdots,j_{w(m)})\in [m]$ to $\mbf i=(i_0,\cdots,i_{w(\mbf
  n)})$ with $i_l=\alpha_l(j_{f^{-1}(l)})$ for $l\in\Image(f)$, and
$i_l=\alpha_l(0)$ otherwise.\\ 
A polysimplicial set is a functor $\mbf
\Lambda^{\op}\to\Set$. Polysimplicial sets
form a category denoted by $\mbf \Lambda^{\circ}\Set$.\\
One considers $\mbf \Lambda$ as a full subcategory of $\mbf
\Lambda^{\circ}\Set$ by the Yoneda functor. If $\C$ is a polysimplicial set
$\mbf \Lambda/\C$ is the category whose objects are morphisms $[\mbf
n]\to\C$ in $\mbf \Lambda^{\circ}\Set$ and morphisms from $[\mbf n]\to \C$
to $[\mbf m]\to \C$ are morphisms $[\mbf n]\to [\mbf m]$ that make the
triangle commute.\\
A polysimplex $x$ of a polysimplicial set $\C$ is said to be \emph{degenerate} if
there is a non isomorphic surjective map $f$ of $\mbf \Lambda$ such that $x$ is the image
by $f$ of a polysimplex of $\C$. Let $\C_{\mbf n}^{\nd}$ be the subset of
non degenerate polysimplices of $\C_{\mbf n}$.\\
Thanks to an analog of Eilenberg-Zilber lemma for polysimplicial sets
(\cite[lem. 3.2]{berk2}), a morphism $\C'\to \C$ is bijective if and only
if it maps non degenerate polysimplices to nondegenerate polysimplices and
$(\C')_{\mbf n}^{\nd}\to\C^{\nd}_{\mbf n}$ is bijective for any $\mbf n$.\\
There is a functor $O:\mbf \Lambda^{\circ}\Set\to\Poset$ where $O(\C)$ is
the partially ordered set associated to $\Ob(\mbf\Lambda/\C)$ endowed with
the preorder where $x\leq y$ if there is a morphism $x\to y$ in $\mbf
\Lambda/\C$. As a set, $O(\C)$ coincides with the set of equivalence
classes of nondegenerate polysimplices.\\
A polysimplicial set $\C$ is said \emph{interiorly free} if $\Aut(\mbf n)$
acts freely on $\C^{\nd}_{\mbf n}$. If $\C_1\to\C_2$ is a morphism of
polysimplicial sets mapping nondegenerate polysimplices to nondegenerate
polysimplices such that $O(\C_1)\to O(\C_2)$ is an isomorphism and $\C_2$
is interiorly free, then $C_1\to\C_2$ is an isomorphism.\\
Berkovich also defines a \emph{strictly polysimplicial category} $\Lambda$ whose
objects are those of $\mbf \Lambda$, but with only injective morphisms
between them. The functor $\Lambda\to\mbf \Lambda
\to\mbf \Lambda^{\circ}\Set$ extends to a functor
$\Lambda^{\circ}\Ens\to\mbf \Lambda^{\circ}\Ens$ which commutes with direct
limits (the objects of $\mbf\Lambda^{\circ}\Ens$ will be called \emph{strictly
polysimplicial sets}).\\
Berkovich then considers a functor $\Sigma:\mbf \Lambda\to\mcal Ke$ to the
category of Kelley spaces, \emph{i.e.} topological spaces $X$ such that a
subset of $X$ is closed whenever its intersection with any compact subset
of $X$ is closed. This functor takes $[\mbf n]$ to $\Sigma_{\mbf
  n}=\{(u_{il})_{0\leq i\leq p,0\leq l\leq n_i}\in [0,1]^{[\mbf n]}|\sum_l
u_{il}=1\}$, and takes a map $\gamma$ associated to $(J,f,\alpha)$ to
$\Sigma(\gamma)$ that maps $\mbf u=(u_{jk})$ to $\mbf u'=(u'_{il})$ defined
as follows: if $[\mbf m]\neq [0]$ and $i\notin\Image(f)$ or $[\mbf m]=[0]$
then $u'_{il}=1$ for $l=\alpha_i(0)$ and $u'_{il}=0$ otherwise; if $[\mbf
m]\neq [0]$ and $i\in\Image(f)$, then 
$u'_{il}=u_{f^{-1}(i),\alpha_i^{-1}(l)}$ for $l\in\Image(\alpha_i)$ and
$u'_{il}=0$ otherwise.\\ 
This induces a functor, the \emph{geometric realization}, $|\ |:\mbf
\Lambda^{\circ}\Set\to\mcal Ke$ (by extending $\Sigma$ in such a way that
it commutes with direct limits).\\
There is also a bifunctor
$\sq:\mbf\Lambda^{\circ}\Set\times\mbf\Lambda^{\circ}\Set\to\Lambda^{\circ}\Set$
which commutes with direct limits and defined by
$[(n_0,\cdots,n_p)]\sq[(n'_0,\cdots,n'_{p'})]=[(n_0,\cdots,n_p,n'_0,\cdots,n'_{p'})]$. Thus
$|\C\sq\C'|=|\C|\times|\C'|$ (where the product on the right is the product
of Kelley spaces).\\

If $X$ is strictly polystable over $k$ and $x\in \Str(X)$, $\Irr(X,x)$ will denote
  the metric space of irreducible components of $X$ passing through $x$ where
  $d(X_1,X_2)=\codim_x(X_1\cap X_2)$. Then there is a tuple $[\mbf n]$
  such that $\Irr(X,x)$ is bijectively isometric to $[\mbf n]$, and if
  $[\mbf m]\to[\mbf n]$ is isometric, there exists a unique
  $y\in\Str(X)$ with $y\leqslant x$ and a unique isometric bijection
  $[\mbf m]\to\Irr(X,y)$ such that
\[\begin{array}{ccc} {[\mbf n]} & \to & \Irr(X,x) \\ \uar & & \uar \\ {[\mbf m]}
  & \to & \Irr(X,y)\end{array}\] commutes.\\
The functor which to $[\mbf n]$ associates the set of couples $(x,\mu)$ where
$x\in\Str(X)$ and $\mu$ is a isometric bijection $[\mbf n]\to\Irr(X,x)$
defines a strict polysimplicial set $C(X)$ (and thus a polysimplicial set
$\C(X)$).\\
There is a functorial isomorphism of partially ordered sets $O(\C(X))\simeq\Str(X)$.
\begin{prop}[{\cite[prop. 3.14]{berk2}}]\label{berk314} One has a functor $\C:\mcal
  Pst^{\sm}\to\Lambda^{\circ}\Ens$, such that $\C(X)$ is as previously
  defined if $X$ is strictly polystable and, for every étale surjective morphism $X'\to X$:
\[\C(X)=\Coker(\C(X'\times_XX')\rightrightarrows \C(X')).\]
\end{prop}
This functor extends to a functor $\C$ for strictly polystable fibrations over $K$
of length $l$.\\
Let us assume we already constructed $\C$ for strictly polystable fibrations of
length $l-1$ such that $O(\C(\underline X))=\Str(X_{l-1})$. Let $\underline X:X_l\to
X_{l-1}\to\cdots\to\Spec k$ be a
strictly polystable fibration, and let $\underline
X_{l-1}:X_{l-1}\to\cdots\to\Spec k$. Then for every $x'\leqslant
x\in\Str(X_{l-1})$, one has:
\begin{lem}[{\cite[cor.6.2]{berk2}}]\label{berk62} There is a canonical cospecialization functor $\C(X_{l,x})\to \C(X_{l,x'})$
and if $x\leqslant x'\leqslant x''$, the functor  $\C(X_{l,x})\to \C(X_{l,x''})$
coincides with the composition $\C(X_{l,x})\to \C(X_{l,x'})\to \C(X_{l,x''})$.\end{lem}
This extends to a functor \[D:(\mbf
\Lambda/(\C(\underline X_{l-1})))^{\op}\to\Str(X_{l-1})\to\mbf\Lambda^{\circ}\Ens.\]
Berkovich then defines a polysimplicial set (where we set
$\C=\C(\underline X_{l-1})$):
\[ \C(\underline X)=\C\sq D = \Coker(\coprod_{N_1(\mbf\Lambda/\C)} \mbf\Lambda [\mbf n_y]\sq
D_x\rightrightarrows\coprod_{N_0(\mbf\Lambda/\C)}\mbf\Lambda [\mbf n_x]
\sq D_x). \]
This construction extends to (non necessarily
strictly) polystable fibrations:
\begin{prop}[{\cite[prop 6.9]{berk2}}]\label{berk69}
There is a functor $\C:\mcal Pst_l^{\tps}\to\mbf\Lambda\Ens$ such
  that:
\begin{enumerate}[(i)]
\item  for every étale surjective morphism of polystable fibrations $X'\to X$: 
\[\C(X)=\Coker(\C(X'\times_XX')\rightrightarrows \C(X')).\]
\item $O(\C(\underline X))\simeq\Str(X).$\end{enumerate}
\end{prop}

Berkovich attaches to a polystable fibration $\underline{\fk X}=(\fk
X_l\to\fk X_{l-1}\to\cdots\to\Spf(O_K))$ a subset of the generic
fiber $\fk X_{l,\eta}$ of $\fk X_l$, the \emph{skeleton} $S(\underline{\fk
  X})$ of $\underline{\fk X}$, which is canonically homeomorphic to $|\C(\fk X_s)|$ (see~\cite[th. 8.2]{berk2}), and such that $\fk X_{l,\eta}$
retracts by a proper strong deformation onto $S(\underline{\fk X})$.\\ 

In fact, when $\underline{\fk X}$ is non degenerate---for example
generically smooth (we will only use the results of Berkovich to such
polystable fibrations)---the skeleton of $\underline{\fk X}$ depends only
on $\fk X_l$  according to~\cite[prop. 4.3.1.(ii)]{berk3}; such a formal
scheme that fits into a polystable fibration will be called
\emph{pluristable}, and we will note $S(\fk X_l)$ this skeleton.\\ 
In this case~\cite[prop. 4.3.1.(ii)]{berk3} gives a description of $S(\fk
X_l)$, which is independant of the retraction. For any $x,y\in\fk
X_{l,\eta}$, we write $x\preceq y$ if for every étale morphism $\fk
X'\to\fk X_l$ and any $x'$ over $x$, there exists $y'$ over $y$ such that
for any $f\in O(\fk X_\eta)$, $|f(x')|\leq |f(y')|$ ($\preceq$ is a partial
order on $\fk X_{l,\eta}$). Then $S(\fk X_l)$ is just the set of maximal
points of $\fk X_{l,\eta}$ for $\preceq$.\\ 

The retraction to $S(\underline{\fk X})$ commutes with étale morphisms:
\begin{thm}[{\cite[th. 8.1]{berk2}}]\label{berk81}
One can construct, for every polystable fibration $\underline{\fk X}=(\fk
X_l\stackrel{f_{l-1}}{\to}\cdots\stackrel{f_1}{\to}\fk X_1\to\Spf(O_K))$, a proper strong deformation retraction $\Phi^l:\fk X_{l,\eta}\times [0,l]\to \fk X_{l,\eta}$ of $\fk
X_{l,\eta}$ onto the skeleton $S(\underline{\fk X})$ of $\underline{\fk X}$ such that:\begin{enumerate}[(i)]
\item $S(\underline{\fk X})=\bigcup_{x\in S(\underline{\fk X}_{l-1})}S(\fk X_{l,x})$ (set-theoretic disjoint union), where $\underline{\fk X}_{l-1}=(\fk X_{l-1}\to\cdots\to\Spf(O_K))$;
\item if $\phi: \fk Y\to \fk X$ is a morphism of fibrations in $\mcal
  Pst f^{\et}_l$, one has $\phi_{l,\eta}(y_t)=\phi_{l,\eta}(y)_t$ for every
  $y\in\fk Y_{l,\eta}$.\end{enumerate}\end{thm}
Let us describe more precisely how the retraction is defined.\\
If $\fk X=\Spf O_K\{P\}/(p_i-z_i)$ where $P$ is isomorphic to
$\oplus_{0\leq i\leq p}\mbf N^{n_i+1}$,
$p_i=(1,\cdots,1)\in\mbf N^{n_i+1}$ and $z_i\in O_K$, let $\fGm$ be the
formal multiplicative group $\Spf O_K\{T,\frac{1}{T}\}$ over $O_K$, let us denote for any
$n$ by $\fGm^{(n)}$ the kernel of the multiplication $\fGm^{n+1}\to\fGm$ and
let $\fk G$ be the formal completion of the identity in
$\prod_i\fGm^{(n_i)}$ (it is a formal group). Then $\fk G$ acts on $\fk X$. $G=\fk
G_\eta$ acts then on $\fk X_\eta$. $G$ has canonical subgroups $G_t$ for
$t\in [0,1]$ defined by the inequalities $|T_{ij}-1|\leq t$ where $T_{ij}$
are the coordinates in $G$, which is a quotient of
$\prod_i\Gm^{n_i+1}$. $G_t$ has a maximal point $g_t$.\\
Then for $x\in X$, $x_t=g_t*x$ defines the strong deformation (where $*$ is
the multiplication defined in~\cite[§ 5.2]{berk}).\\
If $\fk X$ is étale over $\Spf O_K\{P\}/(p_i-z_i)$, the action of $\fk G$
extends in a unique way to an action on $X$, and $x_t$ is still defined by
$g_t*x$. For any $\fk X$ polystable over $O_K$, one has thus defined the strong
deformation locally for the quasi-étale topology of $\fk X_\eta^{\an}$, and Berkovich verifies that it indeed descends
to a strong deformation on $\fk X$.\\
For a polystable fibration $\fk X\to\fk X_{l-1}\to\cdots\to \Spf O_K$,
we first assume $\fk X=\Spf B\to\fk X_{l-1}=\Spf A$ with
$B=A\{P\}/(p_i-a_i)$ (this will be called a \emph{standard} polystable morphism),
one first retracts fiber by fiber on $S(\fk X/\fk X_{l-1})$, which are strictly polystable. The image
obtained can be identified with $S=\{(x,\mbf r_0,\cdots,\mbf r_p)\in \fk
X_{l-1,\eta}, r_{i0}\cdots r_{in_i}=|a_i(x)|\}$, one then has a homotopy
$\Psi:S\times [0,1]\to S$ by $\Psi(x,\mbf r_0,\cdots,\mbf
r_p,t)=(x_t,\psi_{n_0}(\mbf r_0,|a_0(x_t)|),\cdots,\psi_{n_p}(\mbf
r_p,|a_p(x_t)|))$, where $\psi_n$ is some strong deformation of
$[0,1]^{n+1}$ to $(1,\cdots,1)\in [0,1]^{n+1}$ defined by Berkovich (we
will just need that
$\psi_n(r_i,t)_k^\lambda=\psi_n(r_i^\lambda,t^\lambda)_k$ for any
$\lambda\in\mbf R^{*+}$ and any $k\in [[0,n]]$)
, and $x_t$ is defined by the
strong deformation of $\fk X_{l-1,\eta}$.\\
If $\fk X\to\fk X'\to\fk X_{l-1}$ is a geometrically elementary composition
of an etale morphism and a standard polystable morphism, $S(\fk X/\fk
X_{l-1})\to S(\fk X'/\fk X_{l-1})$ is an isomorphism, so that we deform $\fk X'$ fiber by
fiber onto $S(\fk X/\fk X_{l-1})$, then we just do the same retraction
as for $S(\fk X'/\fk X_{l-1})$. For an arbitrary polystable fibration
$X\to\cdots\to O_K$, this
defines the retraction locally for the quasi-étale topology of $\fk X_{\eta}$, and Berkovich verifies that it
descends to a deformation of $X$.\\

Berkovich deduces from (\ref{berk81}.(ii)) the following corollary:
\begin{cor}[{\cite[cor. 8.5]{berk2}}]\label{berk85} Let $K'$ be a finite
  Galois extension of $K$ and let $\underline{\fk X}$ be a polystable fibration over   
  $O_{K'}$ with a normal generic fiber $\fk X_{l,\eta}$. Suppose we are given an action
  of a finite group $G$ on $\underline{\fk X}$ over $O_K$ and a Zariski open dense subset $U$ of $\fk X_{l,\eta}$. 
  Then there is a strong deformation retraction of the Berkovich space
  $G\backslash U$ to a closed subset homeomorphic to $G\backslash |\C(\underline{\fk X})|$.\end{cor}
More precisely, in this corollary, the closed subset in question is the image of $S(\underline{\fk X})$ (which is $G$-equivariant and contained in $U$) by $U\to G\backslash U$.\\

Theorem \ref{berk81} also implies that the skeleton is functorial with respect to pluristable morphisms:
\begin{prop}[{\cite[prop. 4.3.2.(i)]{berk3}}] If $\phi:\fk X\to\fk Y$ is a pluristable morphism between nondegenerate pluristable formal schemes over $O_K$, $\phi_\eta(S(\fk X))\subset S(\fk Y)$.\end{prop}
In fact, more precisely, from the construction of $S$, $S(\fk Y)=\bigcup_{x\in S(\fk X)}S(\fk Y_x)$.

\section{Tempered fundamental group of a polystable log scheme}
In this section we define a tempered fundamental group for a polystable
fibration over a field, endowed with some compatible log structure (we will
call this a polystable log fibration). To define our tempered fundamental
group, we will need a notion of ``topological covering'' of a két covering $Z$ of
our polystable log fibration $X\to\cdots \to k$. To do this we will define for any $Z$ a
polysimplicial set $\C(Z)$ over the polysimplicial set $\C(X)$,
functorially in $Z$. Thus if $Z$ is a finite Galois covering of $X$ with Galois
group $G$, there is an action of $G$ on $\C(Z)$ which defines an extension
of groups:
\[1\to\gtop(|\C(Z)|)\to \Pi_Z\to G\to 1.\]
Our tempered fundamental group will be the projective limits of $\Pi_Z$
when $Z$ runs through pointed Galois coverings of $X$.

\subsection{Polystable log schemes}

Let $S$ be a fs log scheme.
\begin{dfn} A morphism $\phi:Y\to X$ of fs log schemes will be said:
\begin{itemize}
\item \emph{standard nodal} if $X$
has an fs chart $X\to\Spec P$ and $Y$ is isomorphic to $X\times_{\Spec \mbf
  Z[P]}\mbf Z[Q]$ with $Q=(P\oplus u\mbf N\oplus v\mbf N)/(u+v=a)$ with $a\in
P$.
\item a \emph{strictly plurinodal morphism of log schemes} if for every point $y\in Y$, there
exists a Zariski open neighborhood $X'$ of $\phi(y)$ and a Zariski open neighborhood $Y'$ of $y$ in
$Y\times_XX'$ such that $Y'\to X'$ is a composition of strict étale
morphisms and standard nodal morphisms.
\item a \emph{plurinodal morphism of log schemes} if, locally for the étale
  topology of $X$ and $Y$, it is strictly plurinodal.
\item
 a \emph{strictly polystable
  morphism of log schemes} if for every  point $y\in Y$, there
exists a Zariski open neighborhood $X'$ of $\phi(y)$, an fs chart $P\to A$ of the
log structure of $X'$ and a Zariski open neighborhood $Y'$ of $y$ in
$Y\times_XX'$ such that $Y'\to Y$ factors through a strict étale morphism
$Y'\to X'\times_{\mbf Z[P]}\mbf Z[Q]$ where
$Q=(P\oplus\bigoplus_{i=0}^{p}<T_{i0},\cdots,T_{in_i}>)/(T_{i0}+\cdots+T_{in_i}=a_i)$
with $a_i\in P$.
\item a \emph{polystable
  morphism of log schemes} if, locally for the étale topology of $Y$ and $X$, it is a
strict polystable morphism of log schemes.
\end{itemize}
A \emph{polystable log fibration} (resp. \emph{strictly polystable log
  fibration}) $\underline X$ over $S$ of length $l$ is a sequence of
polystable  (resp. strictly polystable) morphism of log schemes
$X_l\to\cdots\to X_1\to X_0=S$.\\
A morphism of  polystable log fibrations of length $l$ $\underline f:\underline
X\to\underline Y$ is given by morphisms $f_i:X_i\to Y_i$ of fs log schemes for
every $i$ such that the obvious diagram commutes.\\
A morphism $\underline f$ of polystable fibrations will be said két
(resp. strict étale) if $f_i$ is két (resp. strict étale) for all $i$.\end{dfn}
A polystable (resp. strictly polystable) morphism of log schemes is
plurinodal (resp. strictly plurinodal).\\
A plurinodal morphism is log smooth and saturated.\\
\begin{rem}
In the definition of strictly
polystable morphisms of log schemes, if one chooses any other chart $P'\to
O(X')$, then there is a Zariski open neighborhood $X''$ of $x$ in $X'$ such
that $Y''=Y'\times_{X'}X''\to X''$ factors through a strict étale morphism $Y''\to X''\times_{\mbf Z[P']}\mbf Z[Q']$ where
$Q'=(P'\oplus\bigoplus_{i=0}^{p}<T_{i0},\cdots,T_{in_i}>)/(T_{i0}+\cdots+T_{in_i}=a'_i)$
with $a'_i\in P'$ (indeed, one can assume that $X'$ is affine, so that
$X'=\Spec(A)$. Let us choose $a'_i$ that has the same image in
$\overline M(X')$ as $a_i$, so that, in $M(X')$, $a'_i=a_iu_i$ with $u_i\in
A^*$. Then one
just changes for every $i$, $T_{i0}$ by $T_{i0}u_i$).
\end{rem}

\begin{lem} Let $\phi:Y\to X$ be a plurinodal (resp. strictly plurinodal,
  resp. polystable, resp. strictly pluristable) morphism of schemes, such that $Y$
  has a log regular log structure $M_X$ and $\phi$ is smooth over
  $X_{\text{tr}}$. Then
  $(Y,O_Y\cap j_*O^*_{Y_{X_{\text{tr}}}})\to (X,M_X)$ is a 
  plurinodal (resp. strictly plurinodal, resp. polystable, resp. strictly
  pluristable) morphism of log schemes.\end{lem}
\dem
Let us prove it for the case of a stricly polystable morphism.\\
One can assume that $X=\Spec(A)$ has a chart $\psi:P\to A$ and that
$Y=B_0\times_X\cdots\times_XB_p$ with $B_i=\Spec
A[T_{i0},\cdots,T_{in_i}]/T_{i0}\cdots T_{in_i}-a_i$ with $a_i\in A$. Since
$\phi$ is smooth over $X_{\text{tr}}$, $a_i$ is
invertible over $X_{\text{tr}}$, thus after multiplying $a_i$ by an
element of $A^*$ (we can do that by also multiplying $T_{i0}$ by this
element), we may assume that $a_i=\psi(b_i)$ for some $b_i\in P$. Thus $Y=X\times_{\mbf Z[P]}\mbf Z[Q]$ where
$Q=(P\oplus\bigoplus_{i=0}^{p}<T_{i0},\cdots,T_{in_i}>)/(T_{i0}+\cdots+T_{in_i}=b_i)$
with $b_i\in P$. If we endow $Y$ with the log structure $M_Y$ associated with
$Q$, $Y\to X$ becomes a strict polystable morphism of log schemes. In
particular $Y$ is log regular~(\cite[th. 8.2]{kato2}). Since, the set of points of $Y$ where $M_Y$
is trivial is $Y_{X_{\text{tr}}}$, $M_Y=O_Y\cap
j_*O^*_{Y_{X_{\text{tr}}}}$ according to~\cite[prop. 2.6]{niziol}. 
\findem

\subsection{Polysimplicial set of a két log scheme over a polystable
  log scheme}
We will construct here the polysimplicial set of a két log scheme $Z$ over a
polystable scheme, with which Berkovich already associated a polysimplicial
set. To do this we will study the stratification of an fs log scheme
defined by $\rk(z)=\rk(\overline M^{\gp}_z)$, which corresponds to
Berkovich stratification for plurinodal schemes, and we will show that
étale locally a két morphism $X\to Y$ induces an isomorphism between the
posets of the strata of $X$ and $Y$. This will enables us to define the
polysimplicial set of $Z$ étale locally. We will then descend it so that
it satisfies the same descent property as in proposition~\ref{berk69}.\\

For a polystable (log) fibration $\underline X:X\to\cdots\to \Spec k$,
Berkovich defines a polysimplicial set $\C(\underline X)$. In this part we
want to generalize this construction to any két log scheme $Z$ over $X$.\\
When $\underline X$ is strictly polystable, $\C(Z)$ will be defined so that
for any stratum $x$ of $X$ with generic point $\tilde x$ and any object $x'$ of
$\mbf \Lambda /\C(X)$ over $x$, the objects of $\mbf \Lambda/\C(Z)$ above
$x'$ will be in natural bijection with the preimage of $\tilde x$ in $Z$.\\
When $\underline X$ is not assumed anymore to be strictly polystable, we
define $\C(Z)$ by étale descent.\\ 

Let $Z$ be an fs log scheme, one gets a stratification on $Z$ by saying
that a point $z$ of $Z$ is of rank $r$ if $\rklog(z)=\rk(M^{\gp}_{\bar z}/\mathcal
O_{\bar z})=r$ (where $\bar z$ is some geometric point over $z$ and where $\rk$ is the rank of an abelian group of finite type).\\
The subset of points of $Z$ such that the rank is $\leqslant r$ is an open
subset of $Z$ (\cite[cor 2.3.5]{ogus}). We thus get a good stratification.\\
The strata of rank $r$ of $Z$ are then the connected components of the
subset of points $z$ of order $r$. This is a partition of $Z$, and a strata
of rank $r$ is open in the closed subset of points $x$ of rank $\geqslant
r$. It is endowed with the reduced subscheme structure of $Z$.\\
The set of strata is partially ordered by $x\leqslant y$ if and only if
$y\subset \bar x$. One denotes by $\Str_x(Z)$ the poset of strata below
$X$.\\

If $f:Z' \to Z$ is a két morphism, then $\rklog(x)=\rklog(f(x))$, so the strata of
$Z'$ are the connected components of the preimages of the strata of $Z$.\\

Let $Z$ be a plurinodal log scheme over some log point
$(k,M_k)$ of characteristic $p$ and of rank $r_0$ and $z$ is a point of
$Z$.\\

$\rklog(z)=r_0+\rk(z)$ where $\rk(z)$ is the codimension of the strata
containing $z$ in $Z$ for the Berkovich stratification of plurinodal
schemes. Thus the strata are the same for this stratification and the
stratification of Berkovich.\\

\begin{lem}
Let $Z\to\Spec k$ be a plurinodal morphism of
log schemes over a log point and $Z'\to Z$ be a két morphism. Then the
strata of $Z'$ (with the reduced scheme structure) are normal and
thus irreducible .\\
\end{lem}
We will often denote abusively in the same way a stratum
and its generic point.\\
\dem
One can prove that étale locally.\\
Let us assume that $Z\to \Spec k$ has an fs chart \[\begin{array}{ccc} Z &
  \to & \Spec P \\ \dar & & \dar \\ \Spec k & \to & \Spec M\end{array}.\]
Then, according to~\cite[prop II.2.3.2]{ogus} the strata of $Z$ correspond
to the connected components of the preimage of the different points of
$\Spec P$  (and thus of the faces of $P$) which are above the minimal point
of $\Spec M$.\\
Let $Q$ be a $(p')$-kummer monoid over $P$. According to~\cite[prop I.1.3.2]{ogus},
$\Spec Q \to \Spec P$ is bijective (to a face $F$ of $P$ corresponds the
saturation $F_Q$ of $F$ in $Q$).\\
Let thus $F$ be a face of $P$ and $\mathfrak p$ be the prime ideal
corresponding to $P\backslash \mathfrak p$ (and $F_Q$ and $\mathfrak p_Q$
be their saturation in $Q$).\\
Then the preimage of the closure of a point corresponding to $F$ in $\Spec
P$ by the map $\Spec k[P]\to \Spec P$ corresponds to the subscheme $\Spec
k[P]/k[\mathfrak p]$ of $\Spec k[P]$ according to~\cite[I.3.2]{ogus} (and
this is the closure of the stratum $\Str F$ of $\Spec k[P]$ corresponding to $F$). But
\[\Spec(k[Q]/k[\mathfrak p_Q]) \to \Spec(k[Q]/(k[\mathfrak p]k[Q]))=\Spec k[Q] \times_{k[P]}\Spec(k[P]/k[\mathfrak p])\]
is just the reduced closed embedding (because $k[\mathfrak p_Q]/(k[\mathfrak
p]k[Q])$ is a nilpotent ideal of $k[Q]/(k[\mathfrak p]k[Q])$). Thus the
preimage of the closure of the stratum corresponding to $F$ is the support
of the closed subscheme $\Spec(k[Q]/k[\mathfrak p_Q])$.\\
Moreover, there are, according to~\cite[I.3.2]{ogus}, canonical
isomorphisms of schemes $k[P]/k[\mathfrak p] \to k[F]$ and
$k[Q]/k[\mathfrak p_Q] \to k[F_Q]$ (but the log structure on $k[F]$ is not
the one induced by $F$) and the following commutative diagram of schemes:
\[\begin{array}{ccc} k[P]/k[\mathfrak p] & \simeq & k[F] \\ \dar & & \dar \\ k[Q]/k[\mathfrak p_Q] & \simeq & k[F_Q] \end{array}\]
where $k[F] \to k[F_Q]$ is induced by the embedding of monoids $F \to F_Q$.\\
The preimage of $F$ by $\Spec k[P]\to \Spec P$ (\emph{i.e.} the stratum
$\Str F$ of $\Spec(k[P])$) corresponds then to the open subset $\Spec k[F^{\gp}]$ of $\Spec k[F] \simeq \Spec k[P]/k[\mathfrak p]$.\\
The following cocartesian square
\[\begin{array}{ccc} k[F] & \to & k[F^{\gp}]\\ \dar & & \dar \\ k[F_Q] & \to & k[F^{\gp}_Q] \end{array}\]
is cocartesian. But $\Spec k[F^{\gp}_Q] \to \Spec k[F^{\gp}]$ is
étale. Thus, $(\Str F\times_{\Spec k[P]} \Spec k[Q])^{\red} \to \Str F$ is étale.\\
By pulling back along $Z\to \Spec P$, we get that the morphism from a
stratum of $Z_Q$ to the corresponding stratum of $Z$ is étale.\findem


\begin{lem} Let $Z'\to Z$ be a két morphism and $Z\to \Spec k$ be a strictly
  plurinodal morphism of log schemes, then $Z$ is quasinormal.\end{lem}
\dem

We will show that the closure of the connected preimages of a stratum $x$
of $Z$ are normal. One can do that étale locally.\\

Let $z$ be a point of $Z'$.\\
One can assume that $Z'$ is connected and $Z'\to\Spec k$ have an fs chart exact at $z$ \[\begin{array}{ccc} Z' &
  \to & \Spec P \\ \dar & & \dar \\ \Spec k & \to & \Spec M\end{array},\]
such that $Z'\to \Spec k\times_{\Spec\mbf Z[M]}\Spec\mbf Z[P]$ is étale.\\
Then the preimages of the different points of $\Spec P$ (\emph{i.e.} primes
of $P$) all map to different strata of $Z$, since $Z$ is strictly
plurinodal.\\
Thus, the preimage of $x$ is either empty or the preimage of $\fk p$ for
some prime $\fk p$ in $\Spec P$.\\

$\Spec k\to \Spec k[M]$ is the closed embedding corresponding to the face
$M^*$ of $M$ (and to the prime ideal $M\backslash M^*$). Thus it can be
described as $\Spec k[M]/k[M\backslash M^*] \to \Spec k[M]$. We thus have
the following commutative diagram:
\[\begin{array}{ccc} \Spec k[P]/k[P(M\backslash M^*)] & \to & \Spec k[P] \\
  \dar & & \dar \\ \Spec k[M]/k[M\backslash M^*] & \to & \Spec
  k[M] \end{array}\]
The strata of $\Spec k[P]/k[P(M\backslash M^*)]$ correspond bijectively to
the prime ideals $\mathfrak p$ of $P$ which contains $P(M\backslash M^*)$
and the closed subscheme which is the closure of $\Spec k[P]/k[\mathfrak
p]\simeq \Spec k[F]$ (where $F$ is the face $P\backslash \mathfrak
p$). $\Spec k[F]$ is connected, and according to~\cite[prop I.3.3.1
(2)]{ogus}, since $F$ is a saturated monoid (because $P$ is a saturated
monoid and $F$ is a face of $P$), $\Spec k[F]$ is normal (thus irreducible,
there is a single stratum above $F \in \Spec P$). The closures of the
strata of $\Spec k\times_{k[M]}\Spec k[P]$ are thus normal too.\findem

One can then follow the proof of~\cite[lem 2.10]{berk2} in the case of $Z$
két over a strictly plurinodal log scheme:
\begin{lem}\label{lemetale} Let $Z'\to Z$ an étale morphism with $Z$
  két
  over $Z_0$ and let $Z_0\to k$ be a stricly plurinodal scheme. Let $z'$ be a
  stratum of $Z'$ of image $z$. Then the map
  $\Str_{z'}(Z')\to \Str_z(Z)$ is an isomorphism of posets.\end{lem}
One can refine~\ref{lemetale} by showing the result for $Z'\to Z$ két:
\begin{lem}\label{lemkum} Let $Z\to Z'$ be a két morphism with $Z'$
  két over $Z''$ and let $Z''\to k$ be a strictly plurinodal
  scheme. Let $z$ be a stratum of $Z$ of image $z'$. Then the
  map $\Str_{z}(Z)\to \Str_{z'}(Z')$ is an isomorphism of posets.\end{lem}
\dem 
It suffices to prove it when $Z'=Z''$, since if one knows the result for
$Z'\to Z''$ and $Z\to Z''$, this implies the result for $Z\to Z'$.\\

Indeed, if $Z'_0=\Spec k\times_{k[M]}\Spec k[P]$ and $Z_0=\Spec
k\times_{k[M]}\Spec k[Q]$ where $P\to Q$ is két, then the poset of strata
of $Z'_0$ (resp.$Z_0$) is isomorphic to the poset of faces of $P$
(resp. $Q$) which maps to the face $M^*$ of $M$. But we already know that
the usual map between the strata of $Q$ and the strata of $P$ is an
isomorphism, hence the isomorphism $\Str Z'_0 \to \Str Z_0$, and a fortiori $\Str_{z_0}(Z_0)\to \Str_{z'_0}(Z'_0)$.\\
In the general case, locally in a Zariski neighborhood of $z$ and of $z'$,
one has the following commutative diagram:
\[\begin{array}{ccccc}Z & \leftarrow & U & \to & Z_0 \\ \dar & & \dar &
  & \dar \\ Z' & = & Z' & \to & Z'_0\end{array}\]
where the horizontal arrows are étale: they satisfy lemma~\ref{lemetale}. By the particular case
studied before, one gets the result.
\findem

Let us consider now a strictly polystable log fibration $\underline
X:X\to X_{l-1}\to\cdots\to s$ where $s$ is an fs log point.\\
If $f:Z \to X$ is két, one has a functor $D_Z=(\mbf
\Lambda/\C(\underline X))^{\op} \to \Str(X_s)\to\Ens$ which associates
to a stratum of $X$ the set of connected components of the preimage of
the stratum (thanks to lemma~\ref{lemkum}). One may thus build a polysimplicial set
$\C_{\underline X}(Z)=\C(\underline X)\sq D_Z$ (we will often write $\C(Z)$
instead of $\C_{\underline X}(Z)$). This polysimplicial set is still
interiorly free. Obviously, $O(\C(Z))$ is functorially
isomorphic to $\Str(Z)$.
\begin{rem}
Let $\C\to \C'$ be a morphism of polysimplicial sets. Let $\alpha:S\to
O(\C)$ (resp. $\alpha':S'\to O(\C')$) be a morphism of posets such that
$S_{\leq x}\stackrel{\simeq}{\to} O(\C_{\leq\alpha{x}})$ (resp. $S'_{\leq
    x}\stackrel{\simeq}{\to} O(\C')_{\leq\alpha'(x)}$ for any $x$). One has
  functors $F:\mbf \Lambda /\C\to O(\C)\to \Set$ (resp. $F':\mbf \Lambda
  /\C'\to O(\C')\to \Set$), which defines a polysimplicial set $D=F\sq\C$ (resp.
$D'=F'\sq\C'$).\\
Then any morphism of posets $f:S\to S'$ such that
\[\begin{array}{ccc} S & \to & S' \\ \dar & & \dar \\ O(\C) & \to & O(\C')
\end{array}\]
induces a unique morphism of polysimplicial sets $\underline f:D\to D'$
over $\C\to \C'$ such that $O(\underline f)=f$.\\
Thanks to this, to construct morphisms between the polysimplicial sets of
kummer log schemes over strictly plurinodal log schemes, we will often be
reduced to construct a morphism between the posets of strata.
\end{rem}

If $\underline X'\to\underline X$ is a két morphism of polystable log
fibrations, then $\C_{\underline X}(X'_l)$ is canonically isomorphic to
$\C(\underline X')$.\\
If one has a commutative diagram
\[\begin{array}{ccc} Z & \to & Z' \\ \dar & & \dar \\ \underline X & \to &
  \underline X'\end{array}\]
where $\underline X \to\underline X$ is a két morphism of polystable log
fibration, there is an induced morphism $\C_{\underline
  X}(Z)\to\C_{\underline X'}(Z')$. Since it maps non degenerate
polysimplices to non degenerate polysimplices, if
$\Str(Z)\to\Str(Z')$ is an isomorphism of posets, $\C_{\underline
  X}(Z)\to\C_{\underline X'}(Z')$ is an isomorphism.\\
Let $Z'\to Z$ be a két covering, let $Z''=Z'\times_ZZ'$ and
let $x$ be a stratum of $X_s$, then $D_Z(x)=\Coker(D_{Z''}(x)\rightrightarrows
D_{Z'}(x))$. We deduce from it that
\[\C(Z'')=\Coker(\C(Z')\rightrightarrows \C(Z)).\]

One may also define $\C_{\underline X}(Z)$ for $\underline X$ a general polystable
fibration. Let $\underline X'\to\underline X$ be an étale covering, let $\underline X''=\underline X'\times_{\underline X}\underline
X'$ and let $Z'$ and $Z''$ the pullbacks of $Z$ to $X'$ and
$X''$. then one defines $\C_{\underline X}(Z)=\Coker(\C_{\underline
  X''}(Z'')\rightrightarrows \C_{\underline X'}(Z'))$ (it
does not depend of the choice of $\underline X'$).\\

If $Z'\to Z$ is a surjective két morphism over $\underline X$ and $Z''=Z'\times_ZZ'$,
$\Str(Z)=\Coker(\Str(Z'')\rightrightarrows \Str(Z'))$.\\
One thus gets ($\ket(X)$ denotes the category of két log schemes over $X$):
\begin{prop}\label{ketpolysimpcplx} Let $\underline X$ be a polystable log fibration, one has a
  functor $\C_X:\ket(X)\to(\Lambda)^{\circ}\Ens$ such that:
\begin{itemize}
\item if $Z'\to Z$ is a
  két covering of $\ket(X_s)$,
\[\C(Z)=\Coker(\C(Z'\times_ZZ')\rightrightarrows \C(Z')).\]
\item $O(\C(Z))$ is functorially isomorphic to $\Str(Z)$.\end{itemize}\end{prop}
\begin{rem} If one has a két morphism $\underline Y\to\underline X$ of
  polystable fibrations of length $l$, the polysimplicial complex $\C(Y_l)$ we have just
  define by considering $Y_l$ as két over $X_l$ is canonically isomorphic
  to the polysimplicial complex of the polystable fibration $\C(\underline
  Y)$ defined by Berkovich.
\end{rem}

If $Z$ is quasicompact, then there is a connected két covering $s'\to s$
such that all the strata of $Z_{\tilde s'}$ are geometrically irreducible
and $Z_{s'}\to s'$ is saturated. In particular, for any Kummer morphism of fs log
points $s''\to \tilde{s'}$, $\C(Z_{s''})\to\C(Z_{\tilde s'})$ is an
isomorphism.\\
The polysimplicial complex $\C(Z_{\tilde s'})$ for such an $s$ is denoted by $\Cgeom(Z/s)$.

\subsection{Tempered fundamental group of a polystable log fibration}\label{tfgsp}
Here we define the tempered fundamental group of a log fibration $\underline
X$ over an fs
log point. If $T$ is a két covering of $X$, the topological coverings of
$|\C(T)|$ will play the role of the topoological coverings of $T$.\\
 
Let us start by a categorical definition of tempered fundamental groups
that we will use later in our log geometric situation.\\
Consider a fibered category
$ \mcal D\to \mcal C$ such that:
\begin{itemize}
\item $\mcal C$ is a Galois category,
\item for every
connected object $U$ of $\mcal C$, $\mcal D_U$ is a category equivalent to $\Pi_U\tSet$
for some discrete group $\Pi_U$,
\item if $U$ and $V$ are two objects of $\mcal C$, the functor $\mcal
  D_{U\coprod V}\to \mcal D_U\times\mcal D_V$ is an equivalence,
\item if $f:U\to V$ is a morphism in $\mcal C$, $f^*:\mcal D_V\to \mcal
  D_U$ is exact.
\end{itemize}
Then, one can define a fibered category $\mcal D'\to\mcal C$ such that the
fiber in $U$ is the category of descent data of $\mcal D\to\mcal C$ with
respect to the morphism $U\to e$ (where $e$ is the final element of $\mcal
C$).\\
Assume one has a splitting of $\mcal D/\mcal C$.\\
Let $U$ is a connected Galois object of $\mcal C$ and let $G$ be the Galois
group of $U/e$. Then $\mcal D'_U$ can be described in the following way:
\begin{itemize}
\item its objects are couples
$(S_U,(\psi_g)_{g\in G})$, where $S_U$ is an object of $\mcal D_U$ and $\psi_g:S_U\to
g^*S_U$ is an isomorphism in $\mcal D_U$ such that for any $g,g'\in G$,
$(g^*\psi_g')\circ\psi_g=\psi_{g'g}$ (after identifying $(g'g)^*$ and
$g^*{g'}^*$ by the canonical isomorphism to lighten the notations).
\item a morphism $(S_U,(\psi_g))\to (S'_U,\psi'_g)$ is a morphism $\phi:S_U\to
  S'_U$ in $\mcal D_U$ such that for any $g\in G$,
  $\psi'_g\phi=(g^*\phi)\psi_g$.\end{itemize}
There is a natural functor $F_0: \mcal D'_U\to\mcal D_U$, which maps
$(S_U,(\psi_g))$ to $S_U$.
Let $F_U$ be a fundamental functor $\mcal D_U\to \Set$, such that $\Aut
F_U=\Pi_U$.\\
Let $F=F_UF_0$, and $\Pi'_U=\Aut F$.\\
\begin{prop} the natural functor $\mcal D'_U\to\Pi'_U\tSet$ is an
  equivalence.\end{prop}
\dem
The proof is some kind of categorical analog of the proof of~\cite[th. 1.4.5]{andre1}.\\
Indeed, $\mcal D'_U$ has limits and colimits, every object is the direct
sum of its connected components and morphisms $S\to S'$ correspond
bijectively with connected components $\Gamma$ of $S\times S'$ such that
the induced morphism $\Gamma\to S$ is an isomorphism. Thus the main point
to prove is that for any connected object $S$ of $\mcal D'_U$, $F(S)$ is a
connected $\Pi'_U$-set. Let $\mcal D'_{U,S}$ be the category of $S$-objects of
$\mcal D'_U$.\\
Let $I$ be the set of connected components of $F_0(S)=S_U$, so that
$S_U=\coprod_{i\in I} S_i$, where $S_i$ is a connected object of $\mcal
D_U$. $G$ acts naturally on $I$ and $S$ is connected implies that $I$ is a
connected $G$-set. Then the category $\mcal D'_{U,S}$ is canonically equivalent
to the category $D$ thus defined:
\begin{itemize}
\item objects of $D$ are couples $S=((S_i)_{i\in I},(\psi_{g,i})_{(g,i)\in
    G\times I}$ such that $S_i$ is an object of $\mcal D_{U,S_i}$ and
  $\psi_{g,i}:S_i\to g^*S_{gi}$ is an isomorphism in $\mcal D_{U,S_i}$;
\item morphisms $S\to S'$ are families $(\phi_i)_{i\in I}$ of morphisms
  $S_i\to S'_i$ compatible with $(\psi_{g,i})$.\end{itemize}

$\mcal D_{U,S_i}$ is a category equivalent to $H_i\tSet$ where $H_i$ is a
subgroup of $\Pi_U\tSet$ since $S_i$ is connected.\\
If $x=(i,F_i)$ is a couple where $i\in I$ and $F_i$ is a fundamental functor of
$\mcal D_{U,S_i}$, one defines a functor $F_x\mcal D'_{U,S}\to \Set$ which
maps $S$ to $F_i(S_i)$. If $x,x'$ are two such couples, there exists an
isomorphism between $F_x$ and $F_x'$: let $g\in G$ such that $gi=i'$ and
let $x_0=(i',F_ig^*)$, then $F_ig^*$ and $F'_{i'}$ must be isomorphic
(since two fundamental functors of $\mcal D_{U,S_{i'}}$ are isomorphic) so
that $F_{x'}$ and $F_{x_0}$ are isomorphic too; $S\mapsto F_i(\psi_{g,i})$
is an isomorphism of functors $F_x\to F_{x_0}$.\\
Now, let $s,s'$ be two points in $F(S)$, and let $i,i'$ be there images in
$I$. Then if $S'\to S_i$ (resp. $S'\to S_{i'}$) is an object of $\mcal
D_{U,S_i}$ (resp. $\mcal D_{U,S_{i'}}$), one defines $F_s(S'\to S_i)$
  (resp. $F_{s'}(S'\to S_{i'})$) to be
  the preimage in $F(S')$ of $s$ (resp. of $s'$). Let $x=(i,F_s)$ and
  $x'=(i',F_{s'})$, then there is an isomorphism $\alpha:F_x\to F_{x'}$ of
  functors of $\mcal D'_{U,S}$. But since for any object $T$ of $\mcal
  D'_U$, $F(T)=F_x(T\times S)=F_{x'}(T\times S)$, one gets an automorphism
  $\beta:F\to F$ (defined by $\beta_T=\alpha_{T\times S\to S}$) which is
  easily seen to send $s$ to $s'$.\\
Thus $F(S)$ is a connected $\Pi'_U$-set. \findem

\begin{prop} There is a natural exact sequence \[1\to \Pi_U\to\Pi'_U\to
  G\to 1.\]\end{prop}
\dem
$F_0$ induces a morphism $\Pi_U\to\Pi'_U$.\\
There is also a natural exact functor $F_1:G\tSet\to \mcal D'_U$ which maps a
$G$-set $Y$ to the couple $(Y=\coprod_{y\in Y}\{y\},(\psi_g))$ where $Y$ is
a constant object in $\mcal D_U$ and $\psi_g$ maps $y$ to $g\cdot y$.\\
This functor is fully faithful and $FF_1$ is naturally isomorphic to the
canonical functor $G\tSet\to\Set$, so that it induces a surjection
$\Pi'_U\to G$.\\
One also has a (non exact) functor $H_0:\mcal D'_U\to G\tEns$ which maps
$(S_U,(\psi_g))$ to the set of connected components of the object $S_U$ of
$\mcal D_U$, where the action of $g\in G$ is induced by $\psi_g$. If
$S=(S_U,(\psi_g))$ is a connected object of $\mcal D'_U$ such that
$F_0(S)=S_U$ has a trivial connected component, $S\to F_1H_0(S)$ is an isomorphism, so that $S$ is in the essential
image of $F_1$. Thus the sequence $\Pi_U\to\Pi'_U\to G$ is exact at
$\Pi'_U$.\\
Finally, let $S_{0,U}$ be a connected object of $\mcal D_U$. Let $S_U=\coprod_{g\in G}
g^*S_{0,U}$, and \[\psi_g:S_U=\coprod_{h\in G}
h^*S_{0,U}=\coprod_{hg\in G}
(hg)^*S_{0,U}\stackrel{\simeq}{\to}\coprod_{h\in G}
g^*h^*S_{0,U}=g^*(\coprod_{h\in G}
h^*S_{0,U})=g^*S_{0,U}.\] One thus gets an object $S$ of $\mcal D'_U$ such
that a connected component of $F_0(S)$ is isomorphic to $S_{0,U}$. This
implies that $\Pi_U\to\Pi'_U$ is injective.\findem

If $(U_i,u_i)_{i\in I}$ is a cofinal projective system of pointed Galois objects
(and let $P$ be the corresponding object of pro-$\mcal C$),
one may define $\Btemp(\mcal D/\mcal C,P)$ to be the category $\injLim_i
\mcal D'_{U_i}$. An isomorphism of pro-objects $P\to P'$ induces an
equivalence $\Btemp(\mcal D/\mcal C,P')\to\Btemp(\mcal D/\mcal C,P)$, so
that $\Btemp(\mcal D/\mcal C,P)$ does not depend up to equivalence on the
choice of $(U_i)_i$.
Moreover, if $h\in G_i=\Gal(U_i/e)$ the endofunctor $h^*:\mcal D'_{U_i}\to\mcal
D'_{U_i}$ maps $S=(S_{U_i},\psi_g)$ to
$h^*S=(h^*S_{U_i},\psi_{hg}\psi_h^{-1})$. Then $\psi_h:S_{U_i}\to
h^*S_{U_i}$ defines an isomorphism $S\to h^*S$ functorially in $S$. Thus
$h^*:\mcal D'_{U_i}\to \mcal D'_{U_i}$ is canonically isomorphic to the
identity of $\mcal D'_{U_i}$. Thus every isomorphism of the pro-object $P$
induces an endofunctor of $\Btemp(\mcal D/\mcal C,P)$ which is canonically
isomorphic to the identity (functorially on $\Aut P$).\\

Let $(F_i)_{i\in I}$ be a family of fundamental functors $F_i:\mcal
D_{U_i}\to \Set$ and assume one has a family $(\alpha_f)_{f:U_i\to U_j}$, indexed on the
set of morphisms in $I$, of isomorphisms of functors $F_if^*\to F_j$ such
that for any $U_i\stackrel{f}{\to} U_j\stackrel{g}{\to} U_k$,
$(\alpha_f\cdot g^*)\alpha_g=\alpha_{gf}$ (after identifying $(gf)^*$ and
$f^*g^*$ to lighten the notations). Such a family exists if $I$ is just
$\mbf N$.
Then, this
induces a projective system  $(\Pi'_{U_i})_{i\in I}$ (unique up to isomorphism
independantly of $(\alpha_f)$ if
$I$=$\mbf N$ and the functors $\mcal D'_{U_i}\to\mcal D'_{U_j}$ are fully
faithful), so that one can define
\[\gtemp(\mcal D/\mcal C,(F_i))=\varprojlim\Pi'_{U_i}\]

Assume one has a 2-commutative diagram with fibered vertical arrows:
\[\begin{array}{ccc} \mcal D_1 & \to & \mcal D_2\\ \dar & & \dar\\ \mcal C_1 & \stackrel{f}{\to} &
  \mcal C_2\end{array}\]
such that $f:\mcal C_1\to \mcal C_2$ is exact, and $\mcal D_{1,U}\to\mcal
D_{2,f(U)}$ is exact for every object $U$ of $\mcal C_1$.\\
One then gets a functor $\Btemp(\mcal D_1/\mcal C_1)\to\Btemp(\mcal
D_2/\mcal C_2)$.\\

For example, Let $X$ be a $K$-manifold, $\mcal C$ be the category of
finite étale covering of $X$ and $\mcal D\to \mcal C$ be the fibered
category such that $\mcal D_U$ is the category of topological coverings of
$U$. Then, since finite étale coverings are morphisms of effective descent
for tempered coverings, $\mcal D'_U$ can be identified functorially with the full
subcategory of $\Covtemp(X)$ of tempered coverings $S$ such that $S_U$ is a
topological covering of $U$. If $(U_i,u_i)$ is a cofinal system of pointed
Galois cover of $(X,x)$, then $\Btemp(\mcal C/\mcal D)$ becomes canonically
equivalent with $\Covtemp(X)$.\\

Let us apply our categorical definition of tempered fundamental groups to our log
geometrical case.\\
Let $\underline X:X\to X_{l-1}\to\cdots\to\Spec(k)$ be a polystable
log fibration, and assume that $X$ is connected.\\
Then one has a functor $\Ctop:\KCov(X)\to \Ke$ obtained by
composing the functor $\C$ of proposition~\ref{ketpolysimpcplx} with the
geometric realization functor.\\
One can thus define a fibered category $\Dtop\to \KCov(X)$ such that the
fiber of a két covering of $Y$ of $X$ is the category of topological
coverings of $\Ctop(Y)$ (which is equivalent to
$\gtop(\Ctop(Y))\tEns$).\\
One defines a fibered category $\DDtemp\to\KCov(X)$ such that the
fiber of a két covering $f:Y\to X$ is the category of descent data of
$\Dtop\to \KCov(X)$ with respect to $Y\to X$ (this corresponds
heuristically to the ``tempered'' coverings of X that become topological
after pullback by $Y\to X$).\\
Let $x$ be a log geometric point of $X$ and let $(Y,y)$ be a
log geometrically pointed
connected Galois két covering of $(X,x)$.\\
Let $\tilde y\to |\C(Y)|$ be the closed cell (which
is contractible) of $|\C(Y)|$ that corresponds to the stratum of $Y$
containing $y$. Then one has a fundamental functor $F_y:\Dtop_Y\to\Ens$
that corresponds to the base point $\tilde y$ ($F_y(S)$ is the set of connected
components of $S\times_{|\C(Y)|}\tilde y$). Moreover, for any morphism $f:(Y',y')\to
(Y,y)$, the two functors $F_{y'}f^*$ and $F_y$ are canonically isomorphic.\\
Then one can consider the functor
$F_{(Y,y)}:\DDtemp_Y\to\Ens$ which associates to a descent data $T$ the set
$F_y(T_Y)$.\\
The induced functor $\DDtemp_Y\to \Aut(F_{(Y,y)})\tEns$ is an equivalence
of categories.\\
One has an exact sequence:
\[1\to \gtop(|\C(Y)|,\tilde y)\to \Aut(F_{(Y,y)}) \to \Gal(Y/X)\to 1.\]
Then one defines
\[\gtemp(X,x)^{\mbb L}=\varprojlim_{(Y,y)} \Aut(F_{(Y,y)}),\]where the projective
limit is taken over the directed category $\LGalKCov(X,x)$ of pointed
connected Galois $\mbb L$-finite két
coverings of $(X,x)$.\\

If $x_1\to x_2$ is a specialization of log geometric points of $X$, it induces a natural
equivalence between the category of pointed coverings of $(X,x_1)$ and the
category of pointed coverings of $(X,x_2)$ (we thus identify the two
categories).\\
If $Y$ is a pointed covering $(Y,y_1)$ of $(X,x_1)$, the corresponding
pointed covering of $(X,x_2)$ is $(Y,y_2)$ where $y_2$ is the unique log
geometric point above $x_2$ such that there is a specialization $y_1\to
y_2$ (and this specialization is unique). Then there is a canonical map
$\tilde y_2\to\tilde y_1$ such that
\[\xymatrix{\tilde y_2 \ar[r] \ar[dr] & \tilde y_1 \ar[d]\\ &
  |\C(Y)|}\]commutes.\\
This induces a canonical isomorphism $F_{y_1}\simeq F_{y_2}$, functorial in
$Y$, so that one gets a canonical isomorphism $\gtemp(X,x_1)^{\mbb
  L}\to\gtemp(X,x_2)^{\mbb L}$. If $X$ is connected and $x_1,x_2$ are two
log geometric points of $X$, there exists a sequence of specializations and
cospecializations joining $x_1$ to $x_2$, so that $\gtemp(X,x_1)^{\mbb L}$
and $\gtemp(X,x_1)^{\mbb L}$ are isomorphic.\\

One has an equivalence of categories between
\[\BtempL_{(X,x)}=\injLim \DDtemp_Y
/\LGalKCov(X,x)\] and the category
$\gtemp(X,x)^{\mbb L}\tEns$ of sets with an action of $\gtemp(X,x)^{\mbb L}$ that goes
through a discrete quotient of $\gtemp(X,x)^{\mbb L}$.\\

Assume now that $X$ is log geometrically connected, \emph{i.e.} that
$X_{k'}$ is connected for any két extension $k'$ of $k$.\\
Let $\bar k$ be a log geometric point on $k$, let $\bar x=(\bar x_{k'})$ be
a compatible system of log geometric points of $X_{k'}$ where $k'$ runs
through két extensions of $(k,\bar k)$ (for example, to construct such a
system, one can take a geometric point of $X_{\tilde k}$ where $\tilde k$
is a strict separable closure of $k$; then $\glog(\tilde k)$ is finitely
generated, and one can thus take a countable system of pointed coverings
$\tilde k_i$ of
$\tilde k$, and take geometric points on $X_{\tilde k_i}$ by induction on $i$).\\
Then, one defines $\gtempgeom(X,\bar x)^{\mbb L}=\varprojlim_{k'}
\gtemp(X_{k'},\bar x_{k'})^{\mbb L}$, where
$k'$ runs through két extensions of $k$ in a log geometric point $\bar
k$.\\
Let $\KCovgeom(X)=\injLim \KCov(X_{k'})$ where $k'$ runs through két
extensions of $k$ in $\bar k$. It is the category of log geometric
coverings of $X$.\\
One thus get a fibered category $\Dtopgeom \to \KCovgeom(X)$, whose fiber
in $Y$ is the category of topological coverings of $|\Cgeom(Y)|$.\\
If $Y\to X$ is a log geometric covering, defined over $k'$,
$\Cgeom(Y_{k'})$ does not depend of $k'$, so that one gets a functor
$\KCovgeom(X)\to\Ke$ which maps $Y$ to $|\Cgeom(Y)|$. If $\bar x$ is a
compatible system of points, for any pointed log geometric covering
$(Y,\bar y)$ of $(X,\bar x)$, $\bar y$ defines a fundamental functor
$F_{\bar y}$ of
$\Dtopgeom_Y$ which are canonically isomorphic for any morphism $(Y',\bar
y')\to (Y,\bar y)$.\\
Then $\gtempgeom(X,\bar x)^{\mbb L}=\gtemp(\Dtopgeom/\KCovgeom(X), (F_{\bar y}))^{\mbb L}$.


\section{Comparison result for the pro-$(p')$ tempered fundamental
  group}
If $\underline X:X\to\cdots\to \Spec(O_K)$ is a proper polystable log fibration,
we want to compare the tempered fundamental group of the generic fiber
$X_\eta$ with the tempered fundamental group of the special fiber endowed
with its natural log structure. The specialization theory of the log
fundamental group already gives us a functor from két coverings of the
special fiber and algebraic coverings of the generic fiber. To extend this
to tempered fundamental groups, one has to compare, for any két covering $T_s$
of the special fiber, the topological space $\C(T_s)$ with the Berkovich
space of the corresponding covering $T_{\eta}$ of the generic fiber. Thus
we will define, as in~\cite{berk2}, a strong deformation retraction of
$T_{\eta}^{\an}$ to a subset canonically homeomorphic to $|\C(T_s)|$. We
will construct this retraction étale locally, where $T$ has a Galois
covering $V'$ by some polystable log fibration over a finite tamely
ramified extension of $O_K$. Then the retraction of the tube of $T_s$ is obtained by
descending the retraction of the tube of $V'_s$, defined in
\cite{berk2}. We will then verify that the retraction does not depend on the choice
of $V'$ so that we can descend the retraction we defined étale locally.

\subsection{Skeleton of a két log scheme over a pluristable
  log scheme}\label{retraction}

Let $\underline X:X\to\cdots\to \Spec(O_K)$ be a polystable log fibration over
$\Spec(O_K)$.
\begin{prop}\label{skelretract} For every két morphism $T\to X$, let $\fk
  T_\eta$ be the generic fiber, in the sense of Berkovich, of the
  formal completion of $T$ along its special fiber. Then, there is a
  functorial map
  $|\C(T)_s)|\to \fk T_\eta$, which identifies, $|\C(T_s)|$ with a subset
  $S(T)$ of $\fk T_\eta$ on which $\fk T_\eta$
  retracts by strong deformation.\end{prop}
\begin{rem}
$\fk T_\eta$ is naturally an analytic subdomain of $T_\eta^{\an}$. Moreover
if $T$ is proper over $O_K$ (for example if $X$ is proper, and $T$ is a két
covering), then $\fk T_\eta\to T_\eta^{\an}$ is an isomorphism.
\end{rem}
\dem
Let $f:T\to X$ be a két morphism.\\

For every $x\in T_{s}$, 
let $\underline U:U_l\to\cdots\to U_0$ be a polystable fibration étale over
$\underline X$ such that $(U_l,x_l)$ is an étale neighborhood of $f(x)$,
such that, for every $i$, $U_i$ has an exact chart $P_i\to
A_i$ and compatible morphisms $P_i\to P_{i+1}$ such that the induced
morphism $U_{i+1}\to U_i\times_{\Spec \mbf Z[P_i]}\Spec \mbf Z[P_{i+1}]$ is
strict étale.\\
One has an étale neighborhood
$i:(V,x')\to (T,x)$ of $x$, a $(p')$-Kummer morphism $P_l\to Q$ such
that $V\to X$ factors through a strict étale morphism $V\to U_l\times_{\Spec
  \mbf Z[P_l]}\Spec \mbf Z[Q]$.\\

 Let $P_i \to \frac{1}{n} P_i$ be
the canonical injection. Then, by definition of a $(p')$-Kummer morphism, there exists
$n$ prime to $p$ such that $P_l\to \frac{1}{n}P_l$ factors through $P_l\to Q$. Thus
$V$ has a két Galois covering that comes from a polystable fibration
$\underline U'=V'\to U'_{l-1}\to\cdots\to \Spec O_{K'}$, where $U'_i=U_i\times_{\Spec
  Z[P_i]}\Spec Z[\frac{1}{n}P_]$ for $i\leq l$ and $V'=V\times_{\mbf Z[Q]}\mbf
Z[\frac{1}{n}P]$ (so that there is a strict étale morphism $V'\to U'_l$) over
$O_{K'}$ for some finite tamely ramified extension
$K'=K[\pi^{1/n}]$ of $K$ .
Let us call $G=(\frac{1}{n}P^{\gp}/Q^{\gp})^\vee$ the Galois group of this
két covering.\\

Let us denote by $\fk U, \fk U_i, \fk V, \fk V'$ the formal completions of
$U, U_i, V, V'$ along the special fiber. $\fk V_{\eta}$ will then denote the
generic fiber of $\fk V$ in the sense of Berkovich.\\
The retraction of $\fk V'_{\eta}$ defined in theorem~\ref{berk81} is
$G$-equivariant, so that it
defines a retraction of $\fk V_{\eta}$.\\
Let $S(\ )$ denote the image of the retraction of $(\ )_{\eta}$.
Then $S(\fk V_{\eta})=G\backslash S(\fk V'_{\eta})=G\backslash
|C(V'_s)|=|G\backslash C(V'_s)|=|C(V_s)|$~(corollary~\ref{berk85}).\\

Let us show that the previously defined retraction of $\fk U_{\eta}$ does
not depend on $n$. Let us start by the case of a polystable morphism.\\
Let \[\psi:Z_1=\Spec A[P]/(p_i-\lambda_i)\to Z_2=\Spec A[P]/(p_i-\lambda_i^s) \]
where $P=\mbf N^{|\mbf r|}=\oplus_{(i,j)\in \mbf r}\mbf Ne_{ij}$ and
$p_i=\sum_j e_{ij}$ induced by the multiplication by
$s$ on $P$, where $s$ is an integer prime to $p$ and where $\lambda\in A$.\\
Let $\mbf G$ be the group $\mbf G^{(\mbf r)}_m$ as defined in~\cite[démo du
th. 5.2 étape 2]{berk2}, it acts on $Z_1$ and $Z_2$. One has $\psi(g\cdot x)=g^s\cdot\psi(x)$.\\
Let $T_{ij}$ be the coordinates of $G$. Then $|T^s_{ij}-1|=|T_{ij}-1|$ if $|T_{ij}-1|<1$. Thus, for $t<1$, $(\
)^s:G\to G$ induces an isomorphism $(\ )^s:G_t\to G_t$, and $g^s_t=g_t$.\\
Thus, if $t<1$ (and also for $t=1$ by continuity), \[\psi(x_t)=\psi(g_t\ast x)=g^s_t\ast\psi(x)=g_t\ast\psi(x)=\psi(x)_t.\]
For a standard polystable fibration, the same result easily follows by
induction using that $\psi_n(r_i,t)^{1/s}=\psi_n(r_i^{1/s},t^{1/s})$ (we
kept the notations from the sketch of the proof of theorem~\ref{berk81}).\\
More precisely, suppose we have the diagram:
\[\xymatrix{B=B'[Y_{ij}]/(Y_{i0}\cdots Y_{in_i}-b_i) & B'\ar[l] \\ 
A=A'[X_{ij}]/(X_{i0}\cdots X_{in_i}-a_i) \ar[u]^\phi & A' \ar[u]^{\phi'} \ar[l]
}\]
where $\phi(X_{ij})=Y_{ij}^s$ anf thus $\phi'(a_i)=b_i^s$, and $\tilde\phi':=\Spf\phi':\Spf
B'\to\Spf A'$ is a két morphism of polystable log fibrations and assume by
induction that we already know that $\tilde\phi(x_t)=\tilde\phi(x)_t$.\\
Let $\fk X$ (resp. $\fk X'$, $\fk Y$, $\fk Y'$) denote $\Spf A$
(resp. $\Spf A'$, $\Spf B$, $\Spf B'$).\\
The first part of the retraction of $\fk X_\eta^{\an}$ and $\fk
Y^{\an}_{\eta}$ (consisting of
the retraction fiber by fiber) commutes with $\tilde\phi:=\Spf\phi$
according to the previous case. We thus just have to study the second part
of the retraction.\\
$\tilde\phi$ induces a map:
\[\begin{array}{ccc} S_A & = & \{(x,r_{ij})\in(\fk X')_\eta^{\an}\times [0,1]^{[\mbf
    n]}|r_{i0}\cdots r_{in_i}=|a_i(x)|\}\subset\fk X_\eta^{\an} \\
\dar & & \\
S_B & = & \{(y,r_{ij})\in(\fk Y')_\eta^{\an}\times [0,1]^{[\mbf
    n]}|r_{i0}\cdots r_{in_i}=|b_i(y)|\}\subset\fk Y_\eta^{\an}
\end{array}\]
 which maps $(x,r_{ij})$ to $(\tilde\phi'(x),r_{ij}^{1/s})$ (remark that
 $|a_i(x)|=|b_i(\tilde\phi'(x))|^s$).\\
Then, if $(x,r_{ij})\in S_A$ (we will write $y:=\tilde\phi'(x)$; by
induction assumption, $\tilde\phi'(x_t)=y_t$)
\[\begin{array}{ccl}
\tilde\phi((x,r_{ij})_t) & = &
\tilde\phi((x_t,\psi_{n_i}(r_{ij},|a_i(x_t)|)_k))\\
& = & (y_t,\psi_{n_i}(r_{ij},|a_i(x_t)|)_k^{1/s})\\
& = & (y_t,\psi_{n_i}(r_{ij}^{1/s},|a_i(x_t)|^{1/s})_k)\\
& = & (y_t,\psi_{n_i}(r_{ij}^{1/s},|b_i(y_t)|)_k)\\
& = & (y,r_{ij}^{1/s})_t\\
& = & \tilde\phi(x,r_{ij})_t\end{array}\]
Thus we get that the retraction of $\fk U_{\eta}$ does not depend on $n$.\\

Let $i:W\to T$ be another neighborhood of $x$ satisfying the same
properties as $V$, and $W'$ defined in the same way (one may assume by the
previous remark that we chose the same $n$). Thus,
$W''=V'\times_{T}W'$ is étale over $V'$ and $W'$ (by the canonical projection denoted by $p$ and $p'$).\\
Let thus $y\in \fk V'_{\eta}$ and $y'\in \fk W'_{\eta}$ with same image in
$\fk T_{\eta}$. Let $y''\in \fk W''_{\eta}$ be above $y$ and $y'$.
Then, for every $t\in [0,1]$, $i(y_t)=i(p(\phi(y''))_t)=i(p(\phi(y'')_t))$
by definition of the retraction of $\fk V_{\eta}$. By
using theorem~\ref{berk81}.(ii) again to $\phi$, one gets
$i(y_t)=ip\phi(y''_t)$. By the same argument for $U'$ and since
$ip\phi=i'p'\phi'$, one gets $i(y_t)=i'(y'_t)$. Thus, the retractions of
the different $\fk V_{\eta}$ are compatible and one gets a well defined retraction of
 $\fk T_{\eta}$ (if $(V_i)$ is an étale covering family of $T$, the
 map obtained by glueing the deformation of the different $\fk V_{i,\eta}$
 is continuous since $\coprod \fk V_{i,\eta}\to \fk T_{\eta}$ is
 quasi-étale and surjective and thus a topological factor map).\\

Moreover, if $\phi:T_1\to T_2$
is a két morphism, $\phi(x_t)=\phi(x)_t$. As in theorem~\ref{berk81}.(vi),
it is also compatible with isometric extensions of $K$.\\

If $T$ is covered by $\widetilde V=\bigcup_i V_i$ such that every $V_i$ satisfies the
same property as $V$,
one gets an isomorphism \[S(\fk T_{\eta})=\Coker(S(\fk{\widetilde V\times_T\widetilde
  V}_\eta)\rightrightarrows \fk{\widetilde
  V}))=\Coker(\bigcup_{i,j}|C(V_{i,s}\times_{T_{s}}V_{j,s}|\rightrightarrows
\bigcup_i|C(V_{i,s})|)=|C(T_s)|\]
This isomorphism is functorial in $T$.
\findem

\subsection{Comparison theorem}
Let $K$ be a discrete valuation field. Let $p$ be the residual
characteristic (which can be 0).\\
Let $\underline X:X\to\cdots\to \Spec(O_K)$ a polystable log fibration over
$\Spec(O_K)$.\\
Let us now compare the tempered fundamental group of the
generic fiber, as a $K$-manifold, and the tempered fundamental group of its
special fiber as defined in~§\ref{tfgsp}.\\

\begin{thm}\label{isomfondtemp}
Let $\bar x$ be a geometric point of $X^{\an}_\eta$, and let $\tilde x$ be
its log reduction.
One has a morphism $\gtemp(X^{\an}_\eta,\bar x)^{\mbb L}\to\gtemp(X_s,\tilde
x)^{\mbb
  L}$ which is an isomorphism if $p\notin\mbb L$.\end{thm}
These morphisms are compatible with finite extensions of $K$.\\

\dem
One has two functors $\mbb L\text{-}\KCov(X)\to\mbb
L\text{-}\Covalg(X_{\eta})$, which is an equivalence of categories if $p\notin\mbb
L$, and
$\mbb L\text{-}\KCov(X)\to\mbb L\text{-}\KCov(X_s)$ which is an equivalence
of categories (theorem~\ref{orgsp}). This
gives us an equivalence of categories between $\mbb L\text{-}\KCov(X_s)$ and
$\mbb L\text{-}\Covalg(X_{\eta})$, which enables us to identify the two categories.\\
One has a fibered category $\Dtopan(X)$ over $\mbb L\text{-}\KCov(X)$ whose fiber
at a $\mbb L$-finite két covering $T$ is the category of topological
coverings of $T_\eta^{\an}$. One has also another fibered category
$\Dtops(X)$ over $\mbb L\text{-}\KCov(X)$ obtained by pulling back the
fibered category $\Dtop(X_s)\to\mbb L\text{-}\KCov(X_s)$ defined in
part~\ref{tfgsp} along $\mbb L\text{-}\KCov(X)\to\mbb
L\text{-}\KCov(X_s)$.\\

Proposition~\ref{skelretract} induces an equivalence of fibered categories
$\Dtopan(X)\to \Dtops(X)$.\\

But $\mbb L\text{-}\KCov(X)\to\mbb
L\text{-}\Covalg(X_{\eta})$ induces a morphism \[\gtemp(X_{\eta}^{\an})^{\mbb
  L}\to \gtemp(\Dtopan(X)/\mbb L\text{-}\KCov(X))\] which is an isomorphism
if $p\notin \mbb L$.\\
Similarly, $\mbb L\text{-}\KCov(X)\to\mbb
L\text{-}\KCov(X_s)$ induces an isomorphism \[\gtemp(X_s)^{\mbb
L}\to\gtemp(\Dtops(X)/\mbb L\text{-}\KCov(X))\]
since $\mbb L\text{-}\KCov(X)\to\mbb
L\text{-}\KCov(X_s)$ is an equivalence of categories.\\
\findem

\subsection{Geometric comparison theorem}
We will assume in this section that $p\notin \mbb L$.\\
\begin{thm}\label{isomtempgeom}There is a natural isomorphism
\[\gtempgeom(X_{s})^{\mbb L}\simeq\gtemp(X_{\bar{\eta}})^{\mbb
  L}.\]\end{thm}
\dem
One knows, according to~\cite[prop 5.1.1]{andre2}, that
\[\gtemp(X_{\bar{\eta}}) \simeq \varprojlim_{K_i} \gtemp(X_{K_i}),\]
where $K_i$ runs through the finite extensions of $K$ in $\overline K$.\\
This induces an analog result for the $\mbb L$-version.\\
However, we would like to know, in the case where $p\notin \mbb L$, if one
can only take the projective limit over tamely ramified extensions of $K$
(\emph{i.e.} to két extensions of $O_k$). Then the isomorphism we want
would simply be obtained from theorem~\ref{isomfondtemp} by taking the
porjective limit over két extensions of $O_k$.\\

We have to show that if $T'$ is a $\mbb L$-finite két geometric covering
of $X$ (which is defined over a finite tamely ramified extension of $K$
according to~\cite[prop. 1.15]{kisin}: one can thus assume that $T'$ is
defined over $K$), the universal topological covering $\widetilde T'_\eta$ of $T'_\eta$ is
defined over some tamely ramified extension of $K$.\\
By changing $\Spec O_K$ by some két covering (which
amounts to changing $K$ by some tamely ramified extension) one
may assume that $T'\to\Spec O_K$ is saturated.\\
One already knows that $\widetilde T'_\eta$ is defined over some finite extension
 $K_2$ of $K$ (\cite[lem 5.1.3]{andre2}). Let $K_1$
be the maximal unramified extension of $K$ in
$K_2$. As $T'\to O_K$ is saturated,
the underlying scheme of $T'_{O_{K_2}}$ is obtained by the base change of schemes $\Spec O_{K_2}\to \Spec O_{K_1}$ of
the underlying scheme of $T'_{O_{K_1}}$. By looking at the special fiber, as
$K_1=K_2$ (as schemes), the morphism $T'_{K_2}
\to T'_{K_1}$ induces an isomorphism between the underlying schemes, thus a
bijection between their strata, and thus an isomorphism
$|C(T'_{K_2})|\to |C(T'_{K_1})|$ and $S(T'_{K_2})\to S(T'_{K_1})$.\\
Thus $\widetilde T'_\eta$ is defined over $K_1$.\\ 
\findem

This isomorphism is $\Gal(\bar K,K)$-equivariant (since the isomorphism for
each Galois extension $K_i$ of $K$ is $\Gal(K_i/K)$-equivariant).\\

\begin{rem} We could also have constructed our isomorphism by taking the
  projective limit over all the separable extensions of $K$ by remarking
  that if $K_1\to K_2$ is totally wildly ramified, then the morphism of log
  points $s_2\to s_1$ is kuh and thus $\gtemp(X_{s_2})\to\gtemp(X_{s_1})$
  is an isomorphism.\end{rem} 

\section{Cospecialization of pro-$(p')$ tempered fundamental group}

Let $X\to Y$ be a proper polystable log fibration, such that $Y$ is log
smooth and proper over $O_K$ (the properness of $Y\to O_K$ is only assumed
so that every point of $Y_\eta$ has a reduction in $Y_s$, but the
cospecialization morphisms we will construct only depend of $Y$ locally).
In this section we will construct the cospecialization morphisms for the
$(p')$-tempered fundamental group of the geometric fibers of $X_\eta\to Y_\eta$. Thanks to theorem~\ref{isomtempgeom} we
will be reduced to construct cospecialization morphisms for the
$(p')$-tempered fundamental group of the log geometric fibers of $X_s\to
Y_s$. Let thus $\bar s_2\to\bar s_1$ be a specialization of log
geometric points of $Y$, where $\bar s_1$ and $\bar s_2$ are
the reductions of geometric points $\bar \eta_1,\bar \eta_2 $of
$Y_\eta$.\\

We already have an equivalence of geometric két coverings of $X_{\eta_1}$ and
$X_{\eta_2}$. Now we must compare, for any such két covering  $Z_{s_1}$
corresponding to $Z_{s_2}$ (which extends to $Z_U$), their polysimplicial sets as defined in
proposition \ref{ketpolysimpcplx}. We will assume that $s_2$ is
the generic point of its stratum (if $s_1$ and $s_2$ are in the same
stratum and $\Cgeom(X_{s_2})$ is interiorly free, it will turn out that $\Cgeom(X_{s_1})\to\Cgeom(X_{s_2})$ is an
isomorphism, so that we can replace $s_2$ by the generic point of its
stratum).  We will construct the cospecialization morphism of
polysimplicial set étale locally, so that we can assume $X$ to be strictly
polystable (the properness will not be used for this). This cospecialization morphism of polysimplicial set will be
constructed by associating, after some két localization of the base so that
$Z_U$ becomes saturated, to a stratum $z$ of $Z_{s_1}$ the minimal stratum
$z'$ of $Z_{s_2}$ such that $z$ is in the closure of $z'$ (as in
lemma~\ref{berk62}). Then the closure of $z'$ in the strict localization of the generic point
of $z$ is separable onto its image, so that $z'$ is geometrically
connected, thus defining a geometric stratum of $X_{s_2}$.\\
We will end this article by glueing our specialization isomorphism of $(p')$-log
tempered fundamental group with our cospecialization morphisms of
polysimplicial sets in a cospecialization morphism of tempered fundamental groups.

\subsection{Cospecialization of polysimplicial sets}

In this section, we construct a cospecialization map of polysimplicial set
for a composition of a két morphism and of a log polystable fibration.\\

One can give an analog of~\cite[prop 2.9]{berk2}:
\begin{prop} Let $Z' \to Z$ be a strictly plurinodal morphism of fs log
  schemes, and $Z''\to Z'$ be a két morphism of log schemes. Let
  $z_1$ and $z_2$ two strata of $Z$ (endowed with the inverse image log
  structure $Z$), such that $z_2\leq z_1$ (\emph{i.e.} $z_1 \in \overline{\{ z_2 \} }$). 
Then one has a cospecialization morphism $\Str(Z'_{z_1})\to \Str(Z'_{z_2})$
which maps a stratum $x_1$ of $\Str(Z'_{z_1})$ to the unique maximal
element of $\{ x_2\in \Str(Z'_{z_2}) | x_1\in \bar x_2 \}$.\\
The cospecialization morphism maps minimal points to minimal points.\\
If $z_3\leq z_2\leq z_1$, the obvious diagram of cospecialization
morphism is commutative.\\
\end{prop}
\dem
As for~\cite[prop 2.9]{berk2}, if the result is true for two morphisms
$\phi:Z'' \to Z'$ and $\psi: Z' \to Z$ then it is also true for $\phi
\circ \psi$ because $\Str(Z''_z)=\coprod_{z' \in \Str(Z'_z)}
\Str(Z''_{z'})$. But it is true if $\psi$ is strictly
plurinodal~(\cite[prop 2.9]{berk2}), and it is also true for $\phi$ két according to~(\ref{lemkum}).\findem

\begin{dfn}
One says that a couple of points $(z_2 \leq z_1)$ of a Zariski fs log
scheme $Z$ is a \emph{good couple} if there is some neighborhood $U$ of $z_1$ and
an fs chart $U\to \Spec \mathbf Z[P]$ such that if $\mathfrak p$ is the
image of $z_2$ in $\Spec P$ by $U \to \Spec P$, and if $F=P\backslash
\mathfrak p$, the reduced scheme $\overline{\{z_2\} }$ endowed with the log
structure associated to $F$ (this log scheme will be denoted by
$\overline{\{z_2\} }_F$) by the morphism
\[\overline{\{z_2\}  }\to \Spec k[P]/k[\mathfrak p] \simeq \Spec k[F]\]
is log regular.\\
One says that a couple of geometric points (resp. log geometric points) $(\bar z_2\to \bar
z_1)$ of an fs log scheme $Z$ if there is some
étale (resp. két) neighborhood $U$ of $\bar z_1$ such that $U$ has a global chart
(and thus is Zariski) and the couple $(z_2\leq z_1)$ of underlying
points of $U$ is a good couple.\end{dfn} 
A couple of points $(z_2\leq z_1)$ is a good couple if $z_2$
is the generic point of a stratum of a log regular Zariski fs log scheme
(\cite[prop. 7.2]{kato2}).\\

\begin{lem} Let $Z' \to Z$ be a Zariski log smooth morphism of fs log schemes, let $(z_2,z_1)$ be a
  good couple of points of $Z$. Let $z'_2$ (respectively $z'_1$) be the
  generic point of a stratum of $Z'_{z_2}$ (respectively $Z'_{z_1}$), such
  that $z'_2\leq z'_1$. Then
  $(z'_2,z'_1)$ is a good couple of points of $Z'$. \label{propboncouples}\end{lem}
\dem Since the statement is local on $Z$ and since log regularity is stable under étale base change, one
can assume that one has an fs chart ($\phi: P \to Q$) of $Z\to Z'$ such
that the square
\[\begin{array}{ccc}Z' & \to & \Spec \mathbf Z[Q]\\ \dar & \square & \dar \\ Z & \to & \Spec \mathbf Z[P]\end{array}\]
is cartesian.\\
Let $\mathfrak p$ (respectively $\mathfrak q$) the ideal of $P$
(respectively $Q$) which is the image of $z_1$ (respectively $z_2$) and let
$F=P\backslash \mathfrak p$ (respectively $F'=Q\backslash \mathfrak q$) be
the associated face. We want to show that $\overline{ \{z'_1\} }$ endowed
with the log structure associated to $F'$ is log regular.\\
One has $F=\phi^{-1}(F')$. $\Ker \phi_{|F}^{\gp} \subset \Ker \phi^{\gp}$
and $\Coker \phi_{|F}^{\gp} \subset \Coker \phi^{\gp}$, thus $\phi_{|F}: F
\to F'$ is also a log smooth morphism of monoids.\\ 
Moreover the following diagram of schemes is commutative:
\[\begin{array}{ccccc}\Spec \mathbf Z[F'] & \simeq & \Spec \mathbf Z[Q]/\mathbf Z[\mathfrak q] & \to & \Spec \mathbf Z[Q] \\ \dar & & \dar & & \dar  \\ \Spec \mathbf Z[F] & \simeq & \Spec \mathbf Z[P]/\mathbf Z[\mathfrak p] & \to & \Spec \mathbf Z[P] \end{array}\]
Thus
\[\begin{array}{ccc} \overline{\{ z'_1\} } & \to & \Spec \mathbf Z[F'] \\ \dar & & \dar \\ \overline{\{ z_1\} } & \to & \Spec \mathbf Z[F] \end{array}\]
is commutative. Let $Z''=\overline{\{ z_1\} } \times_{\Spec \mathbf Z[F]}
\Spec \mathbf Z[F']$ and endow it with the log structure associated to
$F'$. Since $Z'' \to \overline{\{ z_1\} }_F$ is log smooth and
$\overline{\{ z_1\} }_F$ is log regular, according
to~\cite[th. 8.2]{kato2}, $Z''$ is log regular.\\
One has a morphism of schemes $\overline{\{ z'_1\} } \to Z''$ which is the
closed embedding of an irreducible component (because $Z''$ is the preimage
of $\overline{\{ \mathfrak q\} }$ in $Z'_{\overline{\{ z_1\} }}$ and $z'_1$
is a generic point of the preimage of $\overline{\{ \mathfrak q\} }$ in
$Z'_{z_1}$ by definition of a stratum), which induces a strict morphism of
log schemes $\overline{\{ z'_1\} }_{F'} \to Z''$. Thus, since $Z''$ is log
regular (and thus normal), $\overline{\{ z'_1\} }_{F'}$ is a
connected component of $Z''$ and thus is log regular too.\findem

\begin{lem}\label{cosp1} Let $\phi: Z' \to Z$ be a strictly plurinodal
  morphism of fs log schemes. Let $(z_2,z_1)$ be a good couple of points of $Z$.\\ 
Then one has a cospecialization morphism $\Str(Z'_{z_1})\to \Str(Z'_{z_2})$
which maps a stratum $x_1$ of $\Str(Z'_{z_1})$ to the unique maximal
element of $\{ x_2\in \Str(Z'_{z_2}) | x_1\in \bar x_2 \}$.\\
The cospecialization morphism maps minimal points to minimal points.\\
If $z_3\leq z_2\leq z_1$, the obvious diagram of cospecialization morphisms
is commutative.\\
\end{lem}
\dem It clearly is true if $\phi$ is a standard morphism $\Spec B \to \Spec
A$ with $f:  P\to A$ a chart of $\Spec A$ and $B=A[u,v]/(uv-f(a))$ where $a
\in M$ (one can for example use~\cite[lem 2.3]{berk2}).\\
By lemma~(\ref{propboncouples}), one sees, as in the proof of~\cite[prop. 2.9]{berk2},
that if the proposition is true for $\phi$ and $\phi'$, it is true for
$\phi \circ \phi'$. Moreover the result is local for the Zariski topology
of $Z'$, thus there is only to show the result for $\phi$ étale, but this
comes from~(\ref{lemetale}) and the fact that $\overline{\{z_2\} }$ is
normal in a neighborhood of $z_1$ and thus the closures in
$Z'\times_Z\overline{\{z_2\}}$ of two different points of $Z_{z_2}$ (which are two
irreducible components of $Z'\times_Z\overline{\{z_2\}}$, which is normal)
have empty intersection. \findem

\begin{rem} If one has a két morphism $Z''\to Z'$ such that $Z''\to Z$ is
  also strictly plurinodal, then, as for~\cite[cor. 2.11]{berk2}, the following diagram
\[\begin{array}{ccc}\Str(Z''_{z_1}) & \to & \Str(Z''_{z_2})\\ \dar & 
  & \dar \\ \Str(Z'_{z_1}) & \to & \Str(Z'_{z_2})\end{array}\]
is commutative.\end{rem}
If $Z'\to Z$ is a strictly polystable morphism, it induces (as
in~\cite[lem. 6.1]{berk2}) a morphism of polysimplicial complexes
\[\C(Z'_{z_1})\to \C(Z'_{z_2}).\]
If now $Z'\to \cdots\to Z$ is a strictly polystable fibration, using
lemma~\ref{propboncouples}, one constructs by induction on the length of the
fibration a morphism of polysimplicial complexes:
\[\C(Z'_{z_1})\to \C(Z'_{z_2}).\]

\begin{rem}The morphism $Z'\to Z$ is saturated, so that if $s_1$ is a log
  geometric point of $Z$ over $\bar z_1$, the morphism
  $\Cgeom(Z'_{\bar z_1}/s_1)\to \C(Z'_{\bar z_1})$ is an
  isomorphism (and the same thing is also true for $z_2$).\end{rem}

Let $(\bar z_2\to\bar z_1)$ be a good couple of log geometric points
of $Z$.\\
By
replacing $Z$ by the strictly local scheme $Z_1$ at $\bar z_1$ (and let us
chose a good chart modeled by $P$ of $Z_1$ at $\bar z_1$), one gets a morphism
$\psi:\Str(Z'_{\bar z_1})\to\Str(Z'_{z'_2})$,
where $z'_2$ is the image of $\bar z_2$ in the strict localization of $Z$
at $\bar z_1$.\\
Let $x_1$ be the generic point of a stratum $\tilde x_1$ of $Z'_{\bar z_1}$, and
let $Z'_1$ be the local scheme of $Z'_{Z_1}$ at $x_1$ and let $Z''_1$ be
the closure of $\psi(x_1)$ in $Z'_1$ ($x_1$ is still in $Z''_1$).\\
However, étale locally in a neighborhood of $x_1$, $Z'\to Z$ is isomorphic to the pullback to $Z$ of $\Spec \mbf Z[Q]\to
\Spec \mbf Z[P]$ where $P\to Q$ is a saturated morphism. Therefore the
morphism from the closure of a stratum $\overline{\{x_2\}}$ of $Z'_{z'_2}$ to its image
$Z_0$ is
étale locally isomorphic to the pullback to $Z_0$ of $\Spec \mbf
Z[F']\to\Spec \mbf Z[F]$ where $F'$ is the face of $Q$ corresponding to
$x_2$ and $F$ is the preimage face of $F'$ in $P$. Then $F\to F'$ is also a
saturated morphism of monoids thanks to lemma~\ref{satface}. Then $\Spec \mbf
Z[F']\to\Spec \mbf Z[F]$ is a separable morphism of schemes.\\
According to~\cite[cor. 18.9.8]{ega4}, the fibers of $Z''_{1}\to
Z_{0}$ are geometrically connected.
In particular the stratum
$x_{2}$, image of $x_1$ by $\psi:\Str(Z'_{\bar
  z_1})\to\Str(Z'_{z'_{2}})$ is  geometrically connected, thus defines a
  stratum of $\Str(Z'_{\bar z_2})$.\\
Thus one gets a canonical morphism $\Str(Z'_{\bar
  z_1})\to\Str(Z'_{\bar z_2})$ which makes the diagram:
\[\xymatrix{ & \Str(Z'_{\bar z_2})\ar[d] \\ \Str(Z'_{\bar z_1})\ar[ur]\ar[r]
  & \Str(Z'_{z'_2})}\]
One thus gets a cospecialization morphism
\[\Cgeom(Z'_{z_1}/z_1)\to\Cgeom(Z'_{ z_2}/z_2).\]

If $Z'\to \cdots\to Z$ is now a polystable fibration of log schemes, after changing $Z$ by
an étale neighborhood, there is $Z''\to Z'$ an étale morphism of polystable
fibrations over $Z$ such that $Z''\to\cdots \to Z$ is a strictly polystable
fibration of log schemes. Then, by taking the cokernel of the horizontal
arrows of the commutative square:
\[\begin{array}{ccc}C(Z'''_{\bar z_1}) & \rightrightarrows &
  \C(Z''_{\bar z_1})\\ \dar & & \dar \\ \C(Z'''_{\bar z_2}) &
  \rightrightarrows & \C(Z''_{\bar z_2})\end{array}\] 
one gets a cospecialization map $\C(Z'_{\bar z_1})\to \C(Z'_{\bar
  z_2})$ (which is functorial with respect to étale morphisms).\\
This cospecialization map commutes with két morphism of polystable log fibrations.

If $Z''\to Z'$ is a két morphism, and $Z'\to Z$ is a polystable
log fibration and let $(\bar z_2\to \bar z_1)$ be a good couple of log geometric
points.
\begin{prop} There is a canonical cospecialization morphism
  \[\Cgeom(Z''_{z_1}/z_1)\to\Cgeom(Z''_{z_2}/z_2)\] which is
  functorial in $Z''$.\\ \end{prop}
\dem 
Let us assume that $Z'\to Z$ is strictly polystable.\\
After replacing $Z$ by some két neighborhood of $\bar z_1$ and
thus changing $z_1$ and $z_2$ by the reduced subscheme of a connected
Galois két covering in $\bar z_1$ and $\bar z_2$, one may assume that $Z$ has a
global chart modeled on a monoid $M$, that there is an étale covering
finite family
$(U'_i)$ of $Z'$, where $U'_i$ is étale over
$Z\times_{\Spec \mbf Z[M]}\Spec \mbf Z[P_i]$, $Z''$ has an étale covering
family $(U''_i)$ where $U''_i$ is étale over $U'_i\times_{\Spec \mbf
  Z[P_i]} \Spec\mbf Z[Q_i]$ and where $M\to Q_i$ is saturated, so that there
is a $(p')$-Kummer morphism of monoid $Q_i\to P'_i$ such that
$U'''_i=U''_i\times_{\Spec\mbf Z[Q_i]}\Spec\mbf Z[P'_i]$ still fits into a
polystable fibration $\underline U'''_i:U'''_i\to\cdots \to Z$ and $\underline
U'''_i\to\underline U'_i$ is a két morphism of log polystable fibrations (in particular $U'''_i\to U''_i$ is a
Galois covering of group $G=((P'_i)^{\gp}/Q_i^{\gp})^{\vee}$, and $Z''\to
Z$ is saturated).\\ 

Then $\C(U'''_{i,z_1})\to\C(U'''_{i,z_2})$ is $G$ equivariant, so
that it induces a morphism
$\C(U''_{i,z_2})=\C(U'''_{i,z_2})/G\to\C(U'''_{i,z_1})/G=\C(U''_{i,z_1})$.
One deduces from the fact that the cospecialization map commutes with két
morphisms of polystable log fibrations that it descends to a morphism
$\C(Z''_{z_1})\to\C(Z''_{z_2})$.\\
By taking the projective limit over strict étale neighborhood of
$\mring{\bar{z_1}}$, one gets a morphism $\Cgeom(Z''_{z_1}/z_1)\to
\C(Z''_{z'_2})$, where $z'_2$ is the image of $\bar z_2$ in the
strict localization of $Z$ at $\mring{\bar{z_1}}$.\\
If one has a log geometrical stratum of $Z''_{z_1}$, using
\cite[cor. 18.9.8]{ega4} as previously, one gets that the image stratum of $Z''_{z_2}$ is geometrically
connected.\\
One thus gets the wanted morphism
$\Cgeom(Z''_{z_1}/z_1)\to\Cgeom(Z''_{z_2}/z_2)$.\\

If $Z'\to Z$ is not assumed to be strictly polystable anymore, after changing $Z$ by
an étale neighborhood, there is an étale morphism $Z'_0\to Z'$ of polystable
log fibrations over $Z$ such that $Z'_0\to\cdots \to Z$ is a strictly polystable
fibration of log schemes. Let $Z'_1=Z'_0\times_{Z'}Z'_0$. Then, the wanted
morphism is obtained by taking the cokernel of the horizontal
arrows of the following commutative square:
\[\begin{array}{ccc}\Cgeom(Z''_{1,z_1}/z_1) & \rightrightarrows &
  \C(Z''_{0,z_1}/z_1)\\ \dar & & \dar \\ \Cgeom(Z''_{1,z_2}/z_2) &
  \rightrightarrows & \Cgeom(Z''_{1,z_2}/z_2)\end{array}\]
\findem

Let us assume now that $Z''\to Z$ is proper, that $Z$ is log regular and
that $\bar z_1$ and $\bar z_2$ are in the same stratum of $Z$ (\ie the
cospecialization map $M_{Z,z_1}\to M_{Z,z_2}$ is an isomorphism). We may
replace $Z$ by its strict localization at $z_1$ (this does not affect our
cospecialization map). In particular $Z$ is log
Zariski. After some further két localization, we will assume $Z''\to Z$ to
be saturated.\\
The fact that $z_1$ and $z_2$ are in the same stratum of $Z$ easily
implies that $\Cgeom(Z''_{z_1})\to\Cgeom(Z''_{z_2})$ maps non degenerate
polysimplices to non degenerate polysimplices (it suffices to look étale locally).\\
Let ${Z''}^{(i)}$ be the closure of $(Z''_{\bar z_2})^{(i)}$ in $Z''$
endowed with the reduced closed subscheme, and
$(\widetilde{Z}'')^{(i)}$ be its normalization. By looking étale locally on
$Z''$ and thanks to the fact that $z_2$ and $z_1$ are in the same stratum,
one sees that $({Z''}^{(i)})_{z_1}$ is just $(Z''_{z_1})^{(i)}$
and that $((\widetilde{Z}'')^{(i)})_{z_1}$ is just the normalization of
$(Z''_{z_1})^{(i)}$. Thus the connected components of $((\widetilde
  Z'')^{(i)})_{z_1}$ and $((\widetilde Z'')^{(i)})_{z_2}$ are in bijections
with the strata of $Z''_{z_1}$ and $Z''_{z_2}$ of rank $i$. As previously
explained, since $Z''\to Z$ is saturated, $(\widetilde Z'')^{(i)}$ is
separable over the closure of $z_2$ and thus the Stein factorizations of
$(\widetilde Z'')^{(i)}$ for all $i$ tells us that $\Str(Z_{\bar
  z_1})\to\Str(Z_{\bar z_2})$ is bijective.\\

If one assumes moreover that $\Cgeom(Z''_{z_2}/z_2)$ is interiorly free
(this is the case if $\Cgeom(Z'_{z_2}/z_2)$ is interiorly free),
\[\Cgeom(Z''_{{z}_1}/z_1)\to\Cgeom(Z''_{z_2}/z_2)\] is also an isomorphism.\\

We still assume $Z''\to Z$ proper, $Z$ log regular, and assume that for any
stratum of $Z$ with geometric generic point $\bar z$, $\Cgeom(Z'_z)$ is
interiorly free.\\
Let $z_2\to z_1$ be any specialization of log geometric points of $Z$. Let
$x\to z_2$ be a specialization where $x$ is a log geometric point over the
generic point of the stratum of $z_2$ (there exists such a specialization). Thus $(x\to z_2)$ and
$(x\to z_1)$ are good couples. Then one has the
morphisms:
\[\Cgeom(Z_{z_2})\stackrel{\simeq}{\leftarrow}\Cgeom(Z_x)\to\Cgeom(Z_{z_1})\]
so that we get a morphism $\Cgeom(Z_{z_2})\to\Cgeom(Z_{z_1})$ which does
not depend on $x\to z_2$ (since any other specialization $x\to z_2$ goes
through our previous specialization).\\ 
Moreover, if $z_1$ and $z_2$ lie in the same stratum, the cospecialization
morphism is an isomorphism.

\begin{rem} Under the same assumptions, let $\bar y_2\to \bar y_1$ be a specialization of log geometric
points over fs log points $y_2\to y_1$ of $Z'$ whose
images $z_2\to z_1$ in $Z$ lie in the same stratum of $Z$.
 One then have a cospecialization
functor
$F:\KCovgeom(Z'_{z_1}/z_1)^{\mbb L}\to\KCovgeom(Z'_{z_2}/z_2)^{\mbb L}$ if
$\mbb L$ does not contain the residual characteristic $p$ at $z_2$ and $z_1$. If
$Z''_{z_1}$ is some geometric két covering of $Z'_{z_1}$, it extends thanks
to corollary~\ref{logsp} to some két neighborhood $U$ of $\bar z_1$
in $Z$. Let $Z''_U\to U$ be this extension (unique after replacing $U$ by
some smaller neighborhood of $\bar z_1$), and $F(Z''_{z_1})$ is just
the fiber in $\bar z_2$ of $Z''_U$. Then one has an isomorphism
$\Cgeom(Z''_1)\simeq\Cgeom(Z''_2)$, which induces a functor of fibered
categories:
\[\begin{array}{ccc}\Dtopgeom(Z'_1) & \to & \Dtopgeom(Z'_2)\\ \dar & &
  \dar\\ \KCovgeom(Z'_{z_1}/z_1)^{\mbb L} & \to & \KCovgeom(Z'_{z_2}/z_2)^{\mbb
    L} \end{array}\]
and thus a specialization morphism $\gtempgeom(Z'_{z_2},y_2)^{\mbb L}
\to\gtempgeom(Z'_{z_1},y_1)^{\mbb L}$.\\
\end{rem}

\subsection{Cospecialization morphisms of pro-$(p')$ tempered fundamental
  groups}

Let $K$ be a discrete valuation field, and $\Spec O_K$ is endowed with its
usual log structure, and assume that the residual characteristic $p$ of $K$ is
not in $\mbb L$ .\\ 
Let $X\to Y$ be a proper and polystable log fibration such that $Y\to \Spec O_K$
is log smooth and proper, and assume that for every geometric point $\bar
y$ of $Y_s$,
$\C(X_{\bar y})$ is interiorly free (this is for example the case if $X\to
Y$ is strictly polystable or if the fibers of $X\to Y$ are semistable curves).\\ 
Let $y_1$ and $y_2$ be two discretely valued points of
$Y_{\triv}^{\an}$ (after replacing $\mcal H(y_1)$ by an isometric
extension, we may assume $\mcal H(y_1)$ to have an algebraically closed
residue field). One has canonical morphisms of fs log schemes $\Spec O_{\mcal
  H(y_i)}\to Y$ for $i=1,2$. Let $s_i$ be the fs log point of $Y$
corresponding to the special point of $\Spec O_{\mcal
  H(y_i)}$ with the inverse image log structure. Let $s'_i$ be the fs log
point of $Y$ with same underlying scheme as $s_i$ but endowed with the
inverse image log structure.\\
A geometric point $\bar y_i$ (of the Berkovich space $Y^{\an}$) over $y_i$ induces a log geometric point
$\bar s_i$ over $s_i$ ($\bar y_i$ can also be seen as a log
geometric point of $Y$ since the log structure of $Y$ is trivial at $y_i$. Let $\bar s'_i$ be a log geometric point of
$s'_i\times_{s_i}\bar s_i$). \\
Let us assume that one has a specialization $\bar s'_2\to\bar
s'_1$ (and there is one as soon as there is a specialization between the
underlying geometric points of $\mring Y$).\\

More precisely, we define $\Pt^{\an}(Y)$ to be the category whose
objects are geometric points $\bar y$ of $Y_{\tr}^{\an}$, such that
$\mcal H(y)$ is discretely valued (where $y$ is the underlying point of
$\bar y$) and $\Hom(\bar y,\bar y')$ is the set of két
specializations from $\bar s$ to $\bar s'$ where $\bar s$
and $\bar s'$ are the log reductions of $\bar y$ and
$\bar y'$.\\
We will also define $\Pt^{\an}_0(Y)$ to be the category obtained from
$\Pt^{\an}(Y)$ by inverting the class of morphisms $\bar
y\to\bar y'$ such that $\bar s$ and $\bar s'$ lie in the
same stratum of $Y$.\\  

\begin{thm} For any morphism $\bar y_2\to\bar y_1$ in
  $\Pt^{\an}(Y)$, there is an outer morphism
\[\gtemp(X_{\bar y_1}^{\an})^{\mbb L}\to\gtemp(X_{\bar
  y_2}^{\an})^{\mbb L},\]
which is an isomorphism if $s_1$ and $s_2$ lie in the same stratum of $Y$.\end{thm}
\dem
One has a cospecialization functor
\[F:\KCovgeom(X_{s_1}/s_1)^{\mbb L}\to\KCovgeom(X_{s_2}/s_2)^{\mbb L}\] which factors through
$\KCovgeom(X_{Z_0}/Z_0)^{\mbb L}$ where $Z_0$ is the strict localization at $s_1$.
Let $\eta$ be some geometric point above the generic point of $Y$.\\
As the cospecialization functor
$\KCovgeom(X_{s_i}/s_i)^{\mbb L}\to\KCovgeom(X_{y_i}/y_i)^{\mbb L}$ and
$\KCovgeom(X_{y_i}/y_i)^{\mbb L}\to\KCovgeom(X_{\eta}/\eta)^{\mbb L}$ are
equivalences (\cite[prop. 1.15]{kisin}), one
gets that $\KCovgeom(X_{s_1}/s_1)^{\mbb L}\to\KCovgeom(X_{s_2}/s_2)^{\mbb L}$ is an
equivalence.\\

If $Z_{s_1}$ is some geometric két covering of $Z_{s_1}$, it extends thanks
to corollary~\ref{logsp} to some két neighborhood $U$ of $\bar s_1$
in $Z$. Let $Z_U\to U$ be this extension (unique after replacing $U$ by
some smaller neighborhood of $\bar s_1$), and $F(Z_{s_1})$ is just
the fiber in $\bar z_2$ of $Z_U$. Then one has a cospecialization morphism
$\Cgeom(Z_{s_1})\to\Cgeom(Z_{s_2})$, which induces a specialization functor
\[\Dtopgeom_{X_{s_2}}(Z_{s_2})\to \Dtopgeom_{X_{s_2}}(Z_{s_1}),\]
which is an equivalence of categories if $\bar s_1$ and $\bar
s_2$ lie in the same stratum of $Y$.\\
Thus we have a 2-commutative diagram:
\[\begin{array}{ccc}\Dtopgeom_{X_{s_1}} & \to & \Dtopgeom_{X_{s_2}}\\ \dar & &
  \dar\\ \KCovgeom(X_{s_1}/s_1)^{\mbb L} & \stackrel{F^{-1}}{\to} &
  \KCovgeom(X_{s_2}/s_2)^{\mbb L} \end{array}\]
where $F^{-1}$ is some quasi inverse of $F$.
This induces a cospecialization morphism
\[\gtempgeom(X_{s_1}/s_1)^{\mbb L}\to\gtempgeom(X_{s_2}/s_2)^{\mbb L}.\]
The comparison morphisms of theorem~\ref{isomtempgeom} gives us the wanted morphism.
\findem

Thus one gets a functor $\Pt^{\an}_0(Y)^{\op}\to \OutGptop$ (where
$\OutGptop$ is the category of topological groups with outer morphisms)
which maps $\bar y$ to $\gtemp(X_{\bar y}^{\an})^{\mbb L}$.\\

Such a functor cannot exist if $p\neq 0$ and $\mbb L$ is the set of all
primes. Indeed, if $X_1$ and $X_2$ are two Mumford curves with isomorphic
stable reduction but with different metrics on the graphs of their stable
models, then their tempered fundamental groups are not isomorphic. Let us
consider a moduli space of stable curves with level structure, endowed with its
canonical log structure, and a geometric point $\bar s$ in the special fiber of the moduli
space such that the corresponding stable curve has totally degenerate
reduction. In particular, it has at least two double points, and thus
the rank of $\overline M_{\bar s}^{\gp}$ is at least two. Let us take two
valuative fs log points $s_1$ and $s_2$ (\emph{i.e.} $\overline
M_{s_i}\simeq\mbf N$) such that the corresponding morphisms $\overline M_{\bar
  s}\to\mbf N$ are linearly independant, for. Let $\eta_1$ and $\eta_2$ be
discretely valued points of the
analytic geometric fiber whose log reductions are $s_1$ and $s_2$. Then the
two corresponding geometric Mumford curves have different metric on the
graph of their stable model, and thus have non isomorphic tempered
fundamental groups. But two geometric log points over $s_1$ and $s_2$ are
isomorphic with respect to specialization for két topology.\\

If one drops the assumption about inner freeness, one still gets this
result, by the same proof:
\begin{thm} For any couple of geometric points $\bar y_1$ and
  $\bar y_2$ of $Y$ above discretely valued points $y_1, y_2$ of
  $Y$. Let $\bar s_1,\bar s_2$ be their log reductions and assume that
  $\bar s_1$ is the generic point of a stratum of $Y$. Then there is a
  functorial cospecialization outer homomorphism:
\[\gtemp(X_{\bar y_1}^{\an})^{\mbb L}\to\gtemp(X_{\bar
  y_2}^{\an})^{\mbb L}.\]\end{thm}

\providecommand{\bysame}{\leavevmode\hbox to3em{\hrulefill}\thinspace}
\providecommand{\MR}{\relax\ifhmode\unskip\space\fi MR }
\providecommand{\MRhref}[2]{%
  \href{http://www.ams.org/mathscinet-getitem?mr=#1}{#2}
}
\providecommand{\href}[2]{#2}

\end{document}